%% file: main.tex
\documentclass{wscpaperproc}
\usepackage{latexsym}
\usepackage{graphicx}
\usepackage{mathptmx}

\usepackage{amsmath}
\usepackage{amsfonts}
\usepackage{amssymb}
\usepackage{amsbsy}
\usepackage{amsthm}
\usepackage{subfigure}
\usepackage{algpseudocode}
\usepackage{algorithm}
\usepackage{algcompatible}
\usepackage{multicol}
\usepackage{wrapfig}
\usepackage{lipsum}
\usepackage{caption}
\usepackage{amsmath}
\DeclareMathOperator*{\argmax}{arg\,max}
\DeclareMathOperator*{\argmin}{arg\,min}

\newcommand{\stitle}[1]{\noindent{\bf #1}}

\usepackage[pdftex,colorlinks=true,urlcolor=blue,citecolor=black,anchorcolor=black,linkcolor=black]{hyperref}

\newtheoremstyle{wsc}
{3pt}
{3pt}
{}
{}
{\bf}
{}
{.5em}
{}

\theoremstyle{wsc}

\begin{document}

%
%

\pagestyle{fancyplain}

\thispagestyle{plain}
\firstPageHead{}

\chead{\fancyplain{}{\itshape Nezami and Anahideh}}

\rhead{}
\cfoot{}
\renewcommand{\headrulewidth}{0pt} 

\input{wscbib.tex}           

\setlength{\baselineskip}{12.7pt}

\title{AN EMPIRICAL REVIEW OF MODEL-BASED ADAPTIVE SAMPLING FOR GLOBAL OPTIMIZATION OF EXPENSIVE BLACK-BOX FUNCTIONS}


\author{Nazanin Nezami\\[12pt]
    Hadis Anahideh\\[12pt]
	Department of Mechanical and Industrial Engineering\\
	University of Illinois at Chicago\\
	842 W Taylor St\\
	Chicago, IL 60607, USA\\}

\maketitle

\section*{ABSTRACT}
This paper reviews the state-of-the-art model-based adaptive sampling approaches for single-objective black-box optimization (BBO). While BBO literature includes various promising sampling techniques, there is still a lack of comprehensive investigations of the existing research across the vast scope of BBO problems.
We first classify BBO problems into two categories: engineering design and algorithm design optimization and discuss their challenges. 
We then critically discuss and analyze the adaptive model-based sampling techniques focusing on key acquisition functions.
We elaborate on the shortcomings of the variance-based sampling techniques for engineering design problems. Moreover, we provide in-depth insights on the impact of the discretization schemes on the performance of acquisition functions. We emphasize the importance of dynamic discretization for distance-based exploration and introduce \emph{EEPA$^{+}$}, an improved variant of a previously proposed Pareto-based sampling technique.
Our empirical analyses reveal the effectiveness of variance-based techniques for algorithm design and distance-based methods for engineering design optimization problems.


\input{introduction}
\input{problems}

\input{preliminaries}
\input{acquisitions}

\input{EEPA+}

\input{experiments}

\input{Result}
\footnotesize

\bibliographystyle{wsc}

\bibliography{demobib}

\section*{AUTHOR BIOGRAPHIES}


\noindent {\bf Nazanin Nezami} is a PhD Student in Mechanical and Industrial Engineering Department at the University of Illinois at Chicago (UIC). She obtained her M.S. degree in Industrial and Systems Engineering from University of Minnesota Twin Cities prior to joining UIC. Her main research interests are Black-Box Optimization, Machine Learning (ML), and Fairness in ML. \\

\noindent {\bf Hadis Anahideh} is a Research Assistant Professor of the Mechanical and Industrial Engineering Department at the University of Illinois at Chicago. She received her Ph.D. degree in Industrial Engineering from the University of Texas at Arlington. Her research objectives center around Black-box Optimization, Sequential Optimization, Active Learning, Statistical Learning, Explainable AI, and Algorithmic Fairness.\\

\newpage

\end{document}

%% file: wscbib.tex
\makeatletter
\let\@internalcite\cite
\def\cite{\def\@citeseppen{-1000}%
    \def\@cite##1##2{(##1\if@tempswa , ##2\fi)}%
    \def\citeauthoryear##1##2##3{##1 ##3}\@internalcite}
\def\citeNP{\def\@citeseppen{-1000}%
    \def\@cite##1##2{##1\if@tempswa , ##2\fi}%
    \def\citeauthoryear##1##2##3{##1 ##3}\@internalcite}
\def\citeN{\def\@citeseppen{-1000}%
    \def\@cite##1##2{##1\if@tempswa, ##2)\else{}\fi}%
    \def\citeauthoryear##1##2##3{##1 (##3)}\@citedata}
\def\citeA{\def\@citeseppen{-1000}%
    \def\@cite##1##2{(##1\if@tempswa , ##2\fi)}%
    \def\citeauthoryear##1##2##3{##1}\@internalcite}
\def\citeANP{\def\@citeseppen{-1000}%
    \def\@cite##1##2{##1\if@tempswa , ##2\fi}%
    \def\citeauthoryear##1##2##3{##1}\@internalcite}
\def\shortcite{\def\@citeseppen{-1000}%
    \def\@cite##1##2{(##1\if@tempswa , ##2\fi)}%
    \def\citeauthoryear##1##2##3{##2 ##3}\@internalcite}
\def\shortciteNP{\def\@citeseppen{-1000}%
    \def\@cite##1##2{##1\if@tempswa , ##2\fi}%
    \def\citeauthoryear##1##2##3{##2 ##3}\@internalcite}
\def\shortciteN{\def\@citeseppen{-1000}%
    \def\@cite##1##2{##1\if@tempswa, ##2\else{}\fi}%
    \def\citeauthoryear##1##2##3{##2 (##3)}\@citedata}
\def\shortciteA{\def\@citeseppen{-1000}%
    \def\@cite##1##2{(##1\if@tempswa , ##2\fi)}%
    \def\citeauthoryear##1##2##3{##2}\@internalcite}
\def\shortciteANP{\def\@citeseppen{-1000}%
    \def\@cite##1##2{##1\if@tempswa , ##2\fi}%
    \def\citeauthoryear##1##2##3{##2}\@internalcite}
\def\citeyear{\def\@citeseppen{-1000}%
    \def\@cite##1##2{(##1\if@tempswa , ##2\fi)}%
    \def\citeauthoryear##1##2##3{##3}\@citedata}
\def\citeyearNP{\def\@citeseppen{-1000}%
    \def\@cite##1##2{##1\if@tempswa , ##2\fi}%
    \def\citeauthoryear##1##2##3{##3}\@citedata}
%
%
%
\def\@citedata{%
    \@ifnextchar [{\@tempswatrue\@citedatax}%
                  {\@tempswafalse\@citedatax[]}%
}

\def\@citedatax[#1]#2{%
\if@filesw\immediate\write\@auxout{\string\citation{#2}}\fi%
  \def\@citea{}\@cite{\@for\@citeb:=#2\do%
    {\@citea\def\@citea{, }\@ifundefined
       {b@\@citeb}{{\bf ?}%
       \@warning{Citation `\@citeb' on page \thepage \space undefined}}%
{\csname b@\@citeb\endcsname}}}{#1}}%

%
\def\@citex[#1]#2{%
\if@filesw\immediate\write\@auxout{\string\citation{#2}}\fi%
  \def\@citea{}\@cite{\@for\@citeb:=#2\do%
    {\@citea\def\@citea{; }\@ifundefined
       {b@\@citeb}{{\bf ?}%
       \@warning{Citation `\@citeb' on page \thepage \space undefined}}%
{\csname b@\@citeb\endcsname}}}{#1}}%

%
\def\@biblabel#1{}
\makeatother



\newdimen\bibindent
\bibindent=0.0em
\def\thebibliography#1{\section*{\refname}\list
   {}{\settowidth\labelwidth{[#1]}
   \leftmargin\parindent
   \itemindent -\parindent
   \listparindent \itemindent
   \itemsep 0pt
   \parsep 0pt}
   \def\newblock{}
   \sloppy
   \sfcode`\.=1000\relax}

%% file: introduction.tex
\section{INTRODUCTION}
\label{sec:intro}

The task of optimizing a set of design parameters of a black-box system is of great importance whether the parameters are for a real engineering process \shortcite{chen2019aerodynamic}, a computer simulation \shortcite{april2003practical}, a real scientific experiment \shortcite{angermueller2020population}, or a complex algorithm \shortcite{feurer2019hyperparameter}. The term black-box is referred to the underlying structure (derivatives or their approximations) of the objective function that cannot be explicitly defined as it is unknown, difficult, or does not exist. 
In most cases, the black-box function is expensive to compute, as the evaluation time and process to obtain the responses of sample points are prohibitive.
The unconstrained black-box optimization (BBO) problem is defined as minimizing or maximizing a black-box function $f(\mathbf{x})$ over a box-constrained input space $\mathbf{x}^L \leq \mathbf{x} \leq \mathbf{x}^U$,
where $\mathbf{x}$ is a $d$-dimensional input vector of variables, $\mathbf{x}^L$ and $\mathbf{x}^U$ are the lower bound and upper-bound of the range of the variables, and $f(\mathbf{x})$ denotes the unknown expensive objective function. Although we use minimization for mathematical derivations, we consider $-f(\mathbf{x})$ for maximization problems.

Derivative-free optimization \shortcite{rios2013derivative} is the most widely used group of techniques in BBO literature. Due to the characteristics of the black-box functions, classical gradient-based optimization methods are not applicable in the BBO context.
Direct-search methods \shortcite{jones1993lipschitzian,alarie2021two} and Model-based methods \shortcite{jiang2020surrogate} are two well-known classes of techniques in derivative-free optimization. The former evaluates function values
in a subset of sample points and selects points based on those values without any derivative approximation. The latter, however, considers a surrogate model to estimate the unknown black-box function. The fitted surrogate is either directly optimized \shortcite{jones2001taxonomy} or is used to guide a sampling core through an adaptive approach to select informative
points in the design space towards the global optimum.

The model-based adaptive sampling strategy is one of the most effective approaches in BBO literature.
Adaptive sampling selects sample points from regions of interest where the \emph{interest} is measured by two key elements of exploration and exploitation. Exploitation attempts to sample points in regions that may
contain global optimum, which can be detected utilizing the properties of the fitted surrogate model in explored regions. 
The exploration, on the other hand, attempts
to find undiscovered regions. 
Designing a sampling core depends on defining acquisition functions to incorporate exploration and exploitation criteria. 
In this paper, we categorize acquisition functions into two main groups of functions that have a single explicit mathematical formula including Expected Improvement (EI) \shortcite{jones1998efficient}, Probability Improvement (PI) \shortcite{kushner1964new}, Entropy Search (ES) \shortcite{hennig2012entropy}, Knowledge Gradient (KG) \shortcite{frazier2008knowledge}, and Upper Confidence Bound (UCB) \shortcite{srinivas2009gaussian}, Dynamic Coordinate Search (DYCORS) \shortcite{regis2007stochastic}, and the multi-criteria Pareto-based including Surrogate Optimization with Pareto Selection (SOP) \shortcite{krityakierne2016sop} and Exploration Exploitation Pareto Approach (EEPA) \shortcite{dickson2014exploration}.
In the former, a control hyperparameter linearly aggregates the exploration and exploitation metrics. In the latter, however, a non-dominated set is obtained from a Pareto front construction on the surrogate estimates as the exploitation dimension and a distance metric (comparing evaluated samples and new candidates) as the exploration metric.
The Pareto construction provides a more flexible trade-off for exploration and exploitation compared to the restricted linear aggregation and does not require tuning of a control parameter.

Furthermore, acquisition functions either use
variance-based exploration (EI, PI, ES, KG, and UCB), or distance-based exploration metrics (SOP, DYCORS, and EEPA).
Although evaluating points with high uncertainty using variance-based acquisition functions is valuable to identify undiscovered regions, it requires accurate
estimates of the error distribution. Thus, it does not provide a global exploration metric. Variance-based methods often converge to local optima because of their limitation in exploring the entire search space effectively using finite samples. As a result, model-independent exploration metrics seem to be more competent for global exploration, either by defining a weighted score acquisition function (e.g., \shortcite{regis2013combining}) or by constructing a Pareto front (e.g, \shortcite{dickson2014exploration}).


In this paper, we aim to investigate the performance of the most commonly used acquisition functions in model-based adaptive sampling approaches for black-box functions on a wide range of applications and conduct a comprehensive comparison to provide useful guidelines for BBO researchers and practitioners. We specifically separate the applications by the nature of the problems as a) \textbf{engineering design} optimization and b) \textbf{algorithmic design} optimization. The first group's ultimate goal is to find a globally optimal design by optimizing the performance objective of an unknown complex system, whereas the second group's goal is to optimize the performance of the learning function.
Our results support the hypothesis that balancing the exploration-exploitation trade-off is more critical for engineering design and global optimization test problems with highly complex and several local optima compared to algorithmic design problems (e.g., hyperparameter tuning) with less oscillated mappings.
The results confirm that the Pareto-based sampling approaches or acquisition functions with explicit emphasis on exploration and exploitation (e.g., UCB) are more suitable techniques for such complex settings. Other variance-based acquisition functions with an uncontrolled emphasis on exploration and exploitation are more successful for less complex settings such as hyperparameter tuning. 

Furthermore, the above-mentioned common acquisition functions are used to sample from a discretized set of generated sample points. Many applications are naturally defined over discrete solution space (e.g., the number of hidden units in Neural networks \shortcite{luong2019bayesian}). For problems with continuous variables ($\mathbf{x} \in \mathbb{R}^n$), the assumption of natural discretized solution space existence is violated.
There are common approaches in the literature to discretize the representation of the continuous solution space. Design of Experiment space-filling methods \shortcite{pronzato2012design} can be used in a one-shot manner to predesign the solution space, or they can be used in a sequential manner to generate discretized candidate points in each iteration. Smarter discretization strategies such as the Dynamic Coordinate Search approach has been proposed by \shortcite{regis2013combining} to generate high fidelity candidate points that are not necessarily space-filling and are rather generated around interesting regions. DYCORS has been shown to be effective in practice compared to one-shot or sequential space-filling approaches. 
EEPA constructs the Pareto front on a randomly generated discretized set of candidates identifying the non-dominated points. Since the performance of EEPA on a fixed predesigned random set is neither practical nor promising, inspired by DYCORS, we design a dynamic discretization scheme for EEPA in this paper, which significantly improves its performance on a variety of problems. We refer to the advanced EEPA as EEPA$^+$.
Based on our observation, these approaches are more sensitive to the discretization scheme.

The results provide deep insights into the choice of acquisition function for the problem of interest for practitioners and shed light on the challenges of each group of techniques for BBO researchers to address them appropriately.

Our empirically justified observations are as follows: 
\vspace{-3mm}
\begin{itemize}
    \item Variance-based are mostly outperformed by distance-based acquisition functions for global optimization test problems. For high-dimensional cases, the outperformance is even more noticeable.  
    \item Distance-based acquisition functions are effective in real-world engineering design problems.
    \item Variance-based acquisition functions are more suitable for hyperparameter tuning problems, although most acquisition functions have close performances. 
    \item Overall, we observed competitive performance for our proposed EEPA$^+$ across all test problems.
\end{itemize}
\vspace{-2mm}

In \S~\ref{sec:bbo}, we elaborate on the applications of BBO and their complexities. In \S~\ref{sec:sampling}, we review the generic model-based BBO and the well-known acquisition functions. \S~\ref{sec:exp} mainly discusses the impact of the discretization scheme for sampling. In \S~\ref{sec:exp}, we perform an extensive experimental evaluation for three main categories of problems utilizing different acquisition functions. \S~\ref{sec:res} and \S~\ref{sec:final} provide discussion on our findings and general guidelines on the performance of each group of techniques on different types of problems.

%% file: problems.tex
\section{BlACK-BOX OPTIMIZATION PROBLEMS}
\label{sec:bbo}

\setlength{\intextsep}{2pt}%
\setlength{\columnsep}{2pt}%
\begin{wrapfigure}{r}{0.4\textwidth}
\vspace{-10mm}
\centering
\subfigure[MNIST-DNN ($d=2$)]{\includegraphics[width=0.36\textwidth]{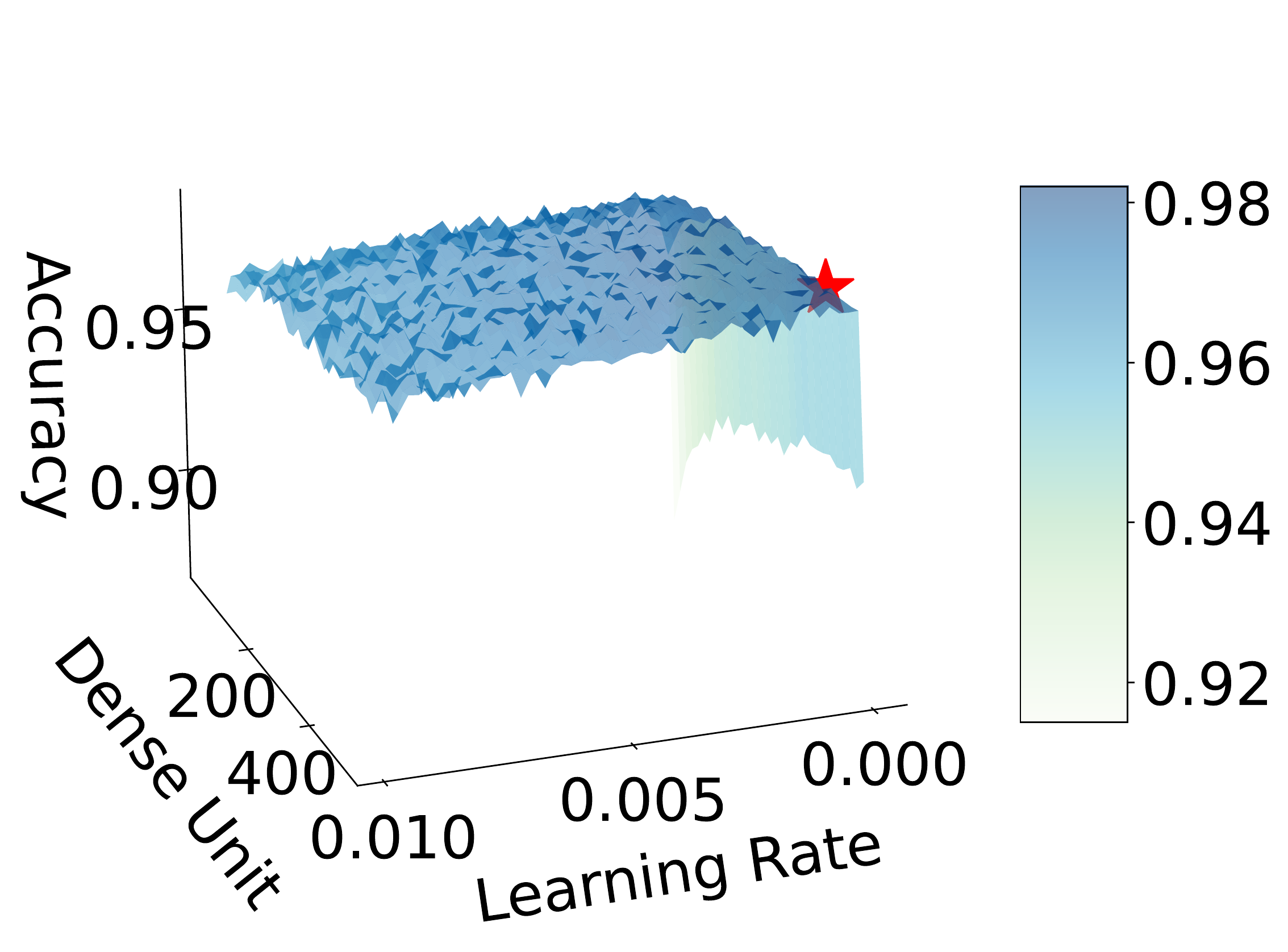},\label{fig:DNN2d}}
\vspace{-3mm}
\subfigure[Rastrigin $d=2$]{\includegraphics[width=0.36\textwidth]{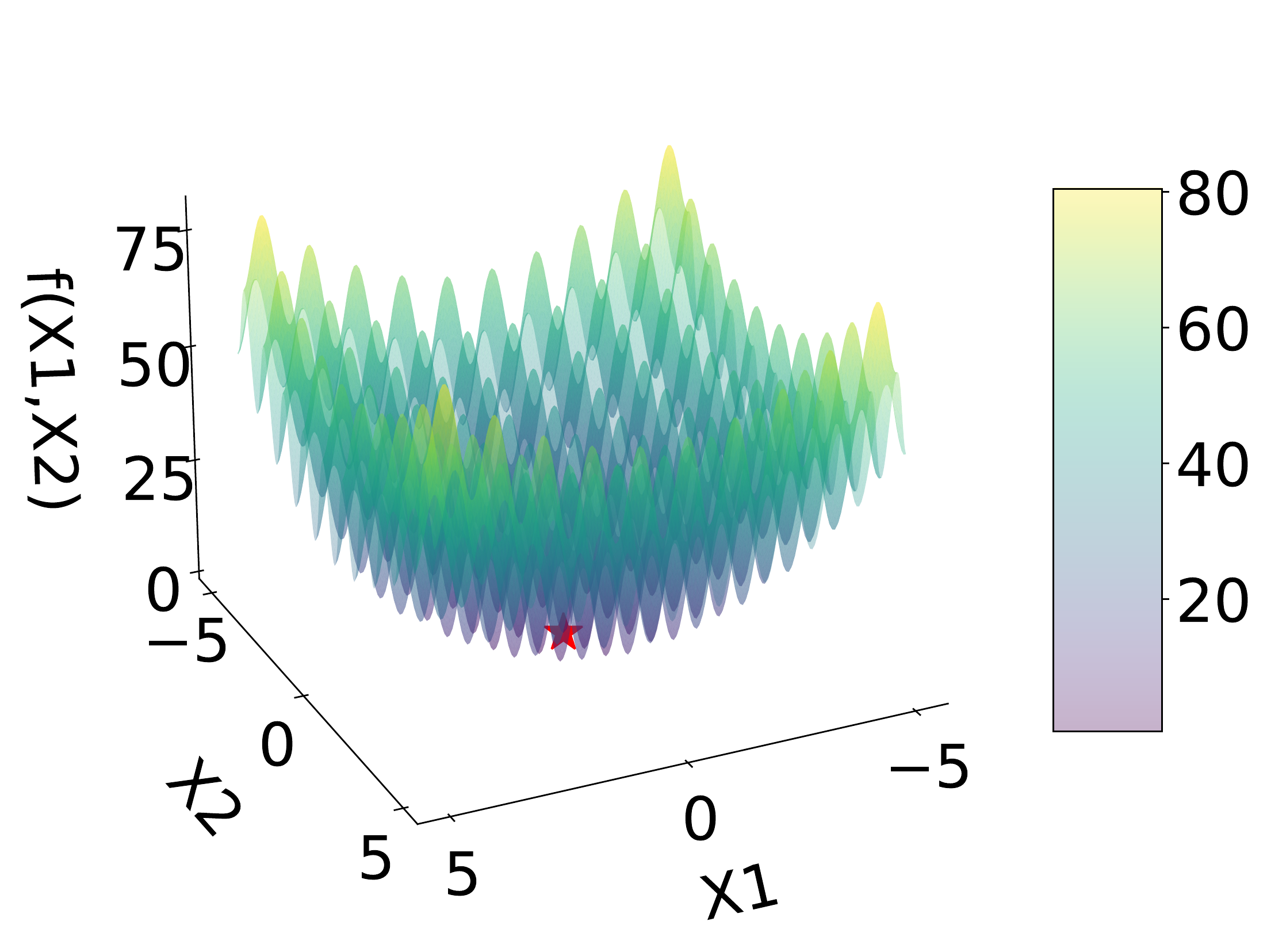},\label{fig:Rastrigin2d}}
\caption{Response surfaces examples}
\vspace{-14mm}
\end{wrapfigure}

Black-box Optimization (BBO) is a well-motivated approach for optimizing expensive to evaluate functions without making specific assumptions about the function behavior. Model-based adaptive sampling strategies for BBO are broadly applicable in many fields with restricted evaluation budgets.
Single objective expensive BBO literature consists of two primary classes of problems; Engineering Design and Algorithm Design optimization. The real-world application of BBO has been investigated in many \textbf{engineering design} problems where the analytical expressions
of the function of interest are inaccessible \shortcite{alarie2021two}. Examples of such problems exist in \emph{airfoil design} \shortcite{palar2019benchmarking,chen2019aerodynamic}, \emph{circuit design} \shortcite{azizi2010energy,lockwood2022empirical}, \emph{groundwater resources} \shortcite{shoemaker2006calibration,muller2021surrogate}, \emph{fluid mechanics} \shortcite{forrester2006optimization,lehnhauser2005numerical}
\emph{chemical processes} \shortcite{griffiths2017constrained,caballero2008rigorous} and other science problems such as \emph{biological design} settings \shortcite{pyzer2018bayesian,angermueller2020population}. 


Another class of BBO problems deals with \textbf{optimizing algorithmic parameters}. These algorithms can be machine learning models or optimization solvers. Optimization solvers (e.g., CPLEX and GUROBI) consist of a large number of combinations of possible solver options that largely affect the performance of the solver \shortcite{audet2006finding,hutter2010automated,liu2019tuning}.
Machine learning models have hyperparameters that must be optimized (aka \emph{Hyperparameter Tuning}) where training time is computationally expensive and the solution space of hyperparameters is often intractable. Instances of BBO for hyperparameter tuning are widespread \shortcite{conn2009introduction,bardenet2013collaborative,wistuba2016hyperparameter,feurer2019hyperparameter}.

Although all of the aforementioned problems are black-box with expensive evaluations, the topology of the unknown underlying response surface may vary around the optimal region. To elaborate on this fact, we provide an example by creating the response surface of the mapping between hyperparameters of the Deep Neural Network (DNN) for the MNIST image classification task through an exhaustive search over a dense grid of two main hyperparameters \emph{(number of dense units, learning rate)} in Figure~\ref{fig:DNN2d}. 
Now, comparing the surface of the DNN tuning with a global optimization test problem (Rastrigin) in Figure~\ref{fig:Rastrigin2d}, on the surface of DNN tuning, we can observe that there are several regions with objective levels comparable to the global optimizer (noted by the red star). This, however, is not the case for the Rastrigin function, which has numerous identifiable local optima and a single unique minimizer. As a result, regret (i.e., $|f(\mathbf{x^*})-f(\mathbf{x})|$) is reduced for hyperparameter tuning due to the modest objective differences between local and global optima, resulting in faster convergence. However, the complex multimodal nature of Rastrigin yields a distinct convergence pattern in which regret minimization can not be easily achieved in areas near global optima. As a result, converging to a global optimum is slower for complex global optimization test problems with identifiable local optima.



%% file: preliminaries.tex
\section{PRELIMINARIES}
\label{sec:pre}

A generic model-based adaptive sampling approach follows an iterative sequential sampling as shown in Algorithm~\ref{alg:Generic}. Let $\mathcal{D}$ be the initial set of evaluated points, commonly obtained from a DOE method, and $\mathcal{F}$ be the set of function values corresponding to the evaluated points. 
The adaptive sampling starts with fitting a cheap to evaluate surrogate model on the evaluated data points. Next, a sample generation scheme (e.g., random) is utilized to discretize the solution space. In the last step, an acquisition function determines the promising candidate points 
with the maximum acquisition value to be evaluated next. Note that $\mathcal{D}$ is updated by appending the newly evaluated points in each iteration. 
The algorithm iterates until the termination criterion is satisfied, which is commonly a maximum evaluation budget.

\textbf{Gaussian processes (GP)} (aka Kriging) \shortcite{rasmussen2003gaussian} is one of the most popular surrogate models in Black-Box Optimization 
that considers the objective function $f(.)$ as the realization of a Gaussian Process.
A Gaussian process (GP), $\hat{f}$, models the predicted mean $\hat{\mu}_t(\mathbf{x})$ and the degree of uncertainty $\hat{\sigma}_t(\mathbf{x})$ at any point $\mathbf{x}$ in the input space, given a set of observations $\mathcal{D}$, where $\mathbf{x}_t$ is the input vector, and $f(\mathbf{x}_t)$ is the corresponding output collected at iteration $t$.
Let $X$ represent the matrix of already evaluated points $\mathcal{D}$, and $\textbf{f}$ be the vector containing elements of $\mathcal{F}$, the estimated mean and variance are updated as follows, 
\vspace{-2mm}
\begin{equation}
\label{eq:GP-mean}
\hat{\mu}(\mathbf{x})=K(X,\mathbf{x} )\Sigma^{-1} \textbf{f}.\\
\end{equation}
\vspace{-5mm}
\begin{equation}
\label{eq:GP-cov}
\hat{\sigma}(\mathbf{x})= K(\mathbf{x} ,\mathbf{x})-K(X,\mathbf{x})\Sigma^{-1}K(\mathbf{x},X)\\
\end{equation}

Where $\Sigma$ represents the covariance matrix and $K$ is a positive semidefinite (PSD) similarity kernel. The off-diagonal elements of $K$ corresponding to any distinct given pair of points $K(\mathbf{x_i},\mathbf{x_j}) ,\forall \mathbf{x_i},\mathbf{x_j} \in \mathcal{X} $ represent the correlation between data points. 
By estimating the prediction errors, GP can guide the sampling core to select points in regions with large variance for reducing the model uncertainty. 

\begin{algorithm}[!hb]
\caption{{\bf Generic Model-Based Adaptive Sampling}}
\begin{algorithmic}[1]
\STATE $\mathcal{D}=\{\mathbf{x}_1, \ldots, \mathbf{x}_{n}\}$,  $\mathcal{F}= \{f(\mathbf{x}_i)|\mathbf{x}_i\in \mathcal{D}\}$ 
\WHILE {Termination criteria not satisfied}
 \STATE \emph{Modeling}: Fit a surrogate model $\hat{f}$ (e.g., GP) on $(\mathcal{D},\mathcal{F})$
 \STATE \emph{Discretization: Generate a set of candidate points, $\mathcal{X}$}
 \STATE \emph{Sampling}: Determine the best candidate points, $S$, using an acquisition function
 \STATE \emph{Evaluation}: $\mathcal{F}_S=\{f(\mathbf{x}_i)|\mathbf{x}_i\in S\}$
  \STATE $\mathcal{D}= \mathcal{D}\cup \mathcal{S}$; $\mathcal{F}=\mathcal{F}\cup \mathcal{F}_S$
\ENDWHILE
\STATE Return $\mathbf{x}^{o}\in\argmin_{\mathbf{x}\in \mathcal{D}} f(\mathbf{x})$ 
\end{algorithmic}\label{alg:Generic}
\end{algorithm}

%% file: acquisitions.tex
\section{ACQUISITION FUNCTIONS FOR MODEL-BASED ADAPTIVE SAMPLING}\label{sec:sampling}
This section briefly discusses the most commonly used acquisition functions for model-based adaptive sampling of black-box optimization.
\subsection{Probability of Improvement (PI)}\label{sec:PI}

Proposed by \shortcite{kushner1964new}, PI is perhaps the first improvement-based acquisition function proposed in black-box optimization literature. Let $\mathcal{D}_n=\{\mathbf{x}_1, \ldots, \mathbf{x}_{n}\}$
represents the set of already evaluated data points and 
$f^{o}_{n}=\max{\{f(\mathbf{x}_i)|\mathbf{x}_i\in \mathcal{D}_n\}}$  be the maximum observed value after $n$ evaluations. The probability of improvement directly measures the probability of inducing improvement over the current best-known solution at a given point $\mathbf{x}$ by defining improvement as $I(\mathbf{x})=\max(f(\mathbf{x})-f^{o}_{n},0)$. Note that based on GP, at any given point, $\mathbf{x}$, the function is sampled from a normal distribution, $f(\mathbf{x})\sim \mathcal{N}(\hat{\mu}(\mathbf{x}),\hat{\sigma}^2(\mathbf{x}))$. Therefore, the probability of improvement can be defined as follows,

\vspace{-5mm}
\begin{equation}
\label{eq:PI}
PI(\mathbf{x})=P(I(\mathbf{x}))=P(f(\mathbf{x})>f^{o}_{n})= \Phi(\frac{\hat{\mu}(\mathbf{x})-f^{o}_{n}}{\hat{\sigma}(\mathbf{x})})
\vspace{-2mm}
\end{equation}
where $\Phi(\mathbf{x})$ indicates the standard normal cumulative distribution function at the point $\mathbf{x}$. 

\subsection{Expected Improvement (EI)}\label{sec:EI}
EI is the most widely-used improvement-based acquisition function in the literature. The EI criterion, which was proposed by Mockus \shortcite{movckus1975bayesian}, gained considerable attention after Jones proposed a closed format expression in \shortcite{jones1998efficient}.
Let $f^{o}_{n}$ and $\mathcal{D}_n$ be the best-known solution and the set of evaluated points after $n$ evaluations, respectively. Now, assuming that for any given point $\mathbf{x}$, the $f(\mathbf{x})|\mathcal{D}_n$ is normally distributed with mean $\hat{\mu}(\mathbf{x})$ and variance $\hat{\sigma}^2(\mathbf{x})$, the closed-form equation of expectation of improvement can be written as, 

\vspace{-5mm}
\begin{equation}
\label{eq:EI}
EI(\mathbf{x})=E(I(\mathbf{x}))= (\hat{\mu}(\mathbf{x})-f^{o}_{n}) \Phi(\frac{\hat{\mu}(\mathbf{x})-f^{o}_{n}}{\hat{\sigma}(\mathbf{x})})+\hat{\sigma}(\mathbf{x})\phi(\frac{\hat{\mu}(\mathbf{x})-f^{o}_{n}}{\hat{\sigma}(\mathbf{x})})
\vspace{-2mm}
\end{equation}
where $I(\mathbf{x})=\max(0,[\hat{\mu}(\mathbf{x})-f^{o}_{n}]^{+})$ defines the EI improvement criteria at a given point $\mathbf{x}$, and $\phi(\mathbf{x})$ and $\Phi(\mathbf{x})$ represent the standard normal density and cumulative distribution functions, respectively. 

\subsection{Upper Confidence Bound (UCB)}\label{sec:UCB}
Upper Confidence Bound (UCB) is one of the most successful 
acquisition functions, which contains explicit objectives for exploration and exploitation. The idea of utilizing UCB criteria for exploration-exploitation trade-off was first mentioned in \shortcite{auer2002using} and the well-known GP-UCB algorithm was later developed by \shortcite{srinivas2009gaussian}.
The exploitation and exploration components are measured by $\hat{\mu}(\mathbf{x})$ and $\hat{\sigma}(\mathbf{x})$ terms, respectively. The UCB equation at any point $\mathbf{x}$ could be described as, 
\vspace{-3mm}
\begin{equation}
\label{eq:UCB}
\alpha_{UCB}(\mathbf{x},\beta_t)=\hat{\mu}(\mathbf{x})+\beta_t^{0.5} \hat{\sigma}(\mathbf{x})
\vspace{-2mm}
\end{equation} 
where $\beta_t$ serves as the trade-off hyperparameter and is commonly set to $2\times log(|\mathcal{X}|t^{2}\pi^{2}/6 \delta)$ formulation which was found to be effective by \shortcite{srinivas2009gaussian}, in which $t$ denotes the number of function evaluations, $|\mathcal{X}|$ represents the cardinality of discretized solution space, and $\delta \in (
0,1)$ is a predefined parameter impacting the probability of achieving the regret bound. Smaller values of $\delta$ increase the probability of achieving the regret bound. 

\subsection{Max-value Entropy Search (MES)}\label{sec:ES}

Entropy Search is among the first information-based acquisition functions that exist in the literature. Proposed by \shortcite{hennig2012entropy}, ES aims at minimizing the uncertainty at the global optimal location $\mathbf{x}^*=\argmax_{x \in \mathcal{X}} f(x)$
by sampling the point that causes the largest entropy reduction. 
Later in \shortcite{hernandez2014predictive}, Predictive Entropy Search (PES) utilized the symmetric property of mutual information and obtained the Entropy Search closed-form expression as $PES(\mathbf{x})=H(P(y|\mathbf{x},\mathcal{D}_n,\mathcal{F}_n))-E[H(P(y|\mathbf{x}^*,\mathcal{D}_n,\mathbf{x})))]$
In other words, while ES describes the expected reduction in entropy for the global optimizer given $\mathbf{x}$, PES measures the expected reduction in entropy for estimated function value $y$ of the new candidate point $\mathbf{x}$. More recently, in \shortcite{wang2017max}, a Max-value Entropy Search (MES) approach has been proposed which is a relatively cheaper and more robust acquisition objective. The key idea is to replace $d$-dimensional $\mathbf{x}^*$ with one-dimensional $y^*$ and rewrite the mutual information equation. MES computes the information gain about the maximum value $y^*=f(\mathbf{x}^*)$ instead of the optimal location $\mathbf{x}^*$ through the following expression:  
\vspace{-2mm}
\begin{equation}
\label{eq:MES}
MES(\mathbf{x})=H(P(y|\mathbf{x},\mathcal{D}_n,\mathcal{F}_n))-E[H(P(y|y^*,\mathcal{D}_n,\mathbf{x})))]
\vspace{-2mm}
\end{equation}


\subsection{Knowledge Gradients (KG)}\label{sec:KG}

Proposed by \shortcite{frazier2008knowledge}, Knowledge Gradients (KG) quantifies the knowledge gain from each candidate point as the conditional expectation of the resulted increase in the posterior mean given that the point is selected for the next evaluation. 
Let $\{(\mathbf{x_i},f(\mathbf{x_i})|i\in{1,\dots,n}\}$ be the set of already evaluated data points, and let $\hat{\mu}_n^{o}=\max_{\forall \mathbf{x} \in \mathcal{X}}{\hat{\mu}_n(\mathbf{x})}$ be the maximum of posterior mean after $n$ evaluations. Now, if we sample any new point, the posterior mean will be updated to $\hat{\mu}_{n+1}^{*}=\max_{\forall \mathbf{x} \in \mathcal{X}}{\hat{\mu}_{n+1}(\mathbf{x})}$ after its evaluation. The knowledge gain from this sampling could then be interpreted as the corresponding increase in the prediction mean for the maximization problem ($\hat{\mu}_{n+1}^{*}-\hat{\mu}_n^{o}$) . As a result, the KG acquisition expression for sampling each point $\mathbf{x}$ for the $n+1$ evaluation could be obtained as the following,
\vspace{-3mm}
\begin{equation}
\label{eq:KG}
KG_n(\mathbf{x})= E_n[\hat{\mu}_{n+1}^{*}-\hat{\mu}_n^{o}|x_{n+1}=\mathbf{x}]
\vspace{-2mm}
\end{equation}



\stitle{Monte Carlo Batch Sampling:}\label{sec:MC} Most analytical acquisition functions, as described in \S~\ref{sec:PI} \S~\ref{sec:EI} \S~\ref{sec:UCB}, require evaluating an expectation term over a discrete finite set of candidate points and returning the best candidate for sampling based on the obtained acquisition values. As a result, they are not directly applicable in a batch setting where joint consideration of sample points is of interest. 
Monte Carlo (MC) sampling approaches can serve as appropriate alternatives to approximate the batch sampling process. 
In this paper, we adopt an efficient MC approach \shortcite{balandat2020botorch} which uses a reparameterized representation to derive approximations for the analytical acquisition functions to sample a batch of sample points. We refer to the associated batched strategies as ``q-acquisition name'' (e.g., qEI). Additionally, we consider the sequential version of these acquisition functions where $q=1$, and denote the batched version of acquisition functions ``s-acquisition name'' (e.g., sEI).  

More specifically, for \textbf{qEI}, \textbf{qPI}, and \textbf{qUCB}, we utilize the reparameterization trick modifying the original 
improvement function $I(.)$ similar to \shortcite{wilson2018maximizing}.
The batch mode of MES acquisition function could be intractable due to the involvement of the multivariate normal cumulative density functions \shortcite{balandat2020botorch}. As a result, solving the optimization in a sequential manner is common practice. 
In our analysis, we utilize a sequential version of \textbf{MES} proposed in \shortcite{wang2017max} which uses Gumbel approximation to sample the max values.  
In addition, for batch KG
we utilize a one-shot approximation, \textbf{qKG}, which converts the problem to a deterministic optimization 
as proposed in \shortcite{balandat2020botorch}).
\subsection{Weighted Score Acquisition Functions}

Originally introduced by \shortcite{regis2007stochastic}, the weighted score acquisition function assigns scores to each candidate point based on a linearly aggregated function over a distance metric and the estimated response value using a surrogate model. Hence, the weighted score considers explicit exploration and exploitation criteria to select the most promising sample points.
Let $w_r^{t}$ and $w_d^{t}$ represent the predefined weight pattern for estimated function value and distance criterion at iteration $t$, respectively. Let $V_r^{t}(\mathbf{x})$ and $V_d^{t}(\mathbf{x})$ be the functions to return a normalized score for the surrogate value and the distance criterion. Let $\Delta(\mathbf{x})_{x \in \mathcal{X}}=\min_{x \in \mathcal{X}}{||\mathbf{x}-x||}$ be the minimum distance from already evaluated data points and let $\hat{f}(\mathbf{x})$ be the GP predicted value. We can obtain the distance score  $V_d^{t}(\mathbf{x})=\frac{\Delta(\mathbf{x})_{max}-\Delta(\mathbf{x})}{\Delta(\mathbf{x})_{max}-\Delta(\mathbf{x})_{min}}$ and the objective score $V_r^{t}(\mathbf{x})=\frac{\hat{f}(\mathbf{x})-\hat{f}(\mathbf{x})_{min}}{\hat{f}(\mathbf{x})_{max}-\hat{f}(\mathbf{x})_{min}}$ in each iteration $t$. As a result, the weighted score acquisition can be written as, 
\vspace{-3mm}
\begin{equation}
\label{eq:score-based}
W^{t}(\mathbf{x})=w_d^{t} \times V_d^{t}(\mathbf{x})+w^{t}_r \times V^{t}_r(\mathbf{x})
\vspace{-2mm}
\end{equation} 
where each $\mathbf{x} \in \mathcal{X}$ represents a candidate point in a 
generated discretized solution space, and $w_d^{t}=1-w_r^{t}$. Note that the above-mentioned description considers the minimization problem by default, however, each minimization setting can be easily converted to the equivalent maximization problem and vice versa.
\vspace{-3mm}
\begin{equation}
\label{eq:problems}
\max f(x) := \min - f(x) 
\vspace{-2mm}
\end{equation}





\subsection{Multi-criteria Pareto-Based Acquisition Function}

\stitle{Exploration Exploitation Pareto Approach (EEPA).}
Proposed by \shortcite{dickson2014exploration}, EEPA is an exploration-exploitation based algorithm that constructs a Pareto frontier over the unevaluated candidate points using two explicit criteria. EEPA utilizes the Maximin distance of the candidate points from the already evaluated points for exploration and the estimated function values for exploitation. More formally $\forall \mathbf{x} \in \mathcal{X}$, the exploration metric for is defined as $\max\limits_{\mathbf{x}\in \mathcal{X}} \min\limits_{\mathbf{x}'\in \mathcal{D}} ||\mathbf{x}-\mathbf{x}'||$, and the exploitation metric is 
denoted as 
$\min\limits_{\mathbf{x} \in \mathcal{X}}{\hat{f}(\mathbf{x})}$. EEPA in \shortcite{dickson2014exploration} utilizes a fixed random set of unevaluated points to represent the solution space. \\

\stitle{Surrogate Optimization with Pareto Selection (SOP).}
SOP, has been introduced by \shortcite{krityakierne2016sop}. Unlike EEPA, SOP constructs a Pareto front based on exploration and exploitation metrics on the already evaluated points, $\mathcal{D}$, rather than unevaluated candidate points, $\mathcal{X}$, to select a few reference points (i.e., \emph{centers}).
The algorithm generates a dynamic discretization of unevaluated sample points through random perturbation of the selected centers. Next, a surrogate-assisted candidate search is performed to select the candidate point with the minimum estimated function value.
A multi-layer heuristic approach is proposed to overcome the shortcomings of single Pareto front construction to generate quality sample points.
SOP creates a Tabu list, which is a list containing the selected centers
that did not significantly improve the hypervolume of the exploration-exploitation trade-off in previous iterations. A point can only be selected as a center if it does not belong to the Tabu list and if its distance from already selected centers is greater than the perturbation radius.

%% file: EEPA+.tex
\section{DISCRETIZATION}
The generalized assumption of continuous input space in BBO is violated in many cases \shortcite{bartz2017model,luong2019bayesian,papalexopoulos2021constrained}. 
Moreover, most of the acquisition functions are only tractable over a finite discretized solution space representation.
Hence, an effective space representation yields a robust search for complex high-dimensional black-box settings. In contrast, a weak representation leads to poor acquisition performance, especially when the black-box is complex and/or higher-dimensional.

Despite the primary role in the model-based adaptive approach, solution space discretization (aka candidate generation) for BBO has only been investigated in a few research articles. 
Using random distributions \shortcite{regis2007stochastic,dickson2014exploration}, a DOE method \shortcite{golzari2015development}, dividing the search space \shortcite{jones1993lipschitzian,eric2007active}, and adaptive grid construction \shortcite{bardenet2010surrogating} to generate a large set of unevaluated points are common practice yet ineffective in most BBO applications. As a result, more practical discretization scheme such as Monte Carlo \shortcite{snoek2012practical} and dynamic \shortcite{regis2013combining,muller2019surrogate} discretization have been proposed in the literature. 
Dynamic discretization methods smartly generate sample points in promising regions and thus, are more effective for high-dimensional complex problems. For example, in \shortcite{regis2013combining}, a dynamic coordinate search strategy, \textbf{DYCORS}, has been proposed, which uses a decreased number of the input coordinates to perturb the current best-known solution and generate promising candidate points. Similarly, \shortcite{muller2019surrogate} utilized a perturbation of the current best-known solution where the perturbation range is drawn from a random normal distribution.
\shortcite{anahideh2022high} introduced a new discretization strategy in which a set of promising candidates, centroids obtained from Decision Tree partitioning, are dynamically added to a randomly generated set of unevaluated points to improve the representation of the solution space with high-quality sample points.
\\

\stitle{EEPA$^+$:An Improved Exploration-Exploitation Pareto Approach.}
Since the original \textbf{EEPA} algorithm \shortcite{dickson2014exploration} utilizes a fixed predesigned discretization of unevaluated sample points, in this paper, we first propose an improved version of \textbf{EEPA} utilizing a dynamic coordinate discretization scheme similar to \textbf{DYCORS} \shortcite{regis2013combining}, which we refer to it as \textbf{EEPA$^+$}. We then empirically investigate the impact of the solution space discretization scheme on the quality of the best-known solution obtained from different acquisition functions.

\begin{algorithm}[!htb]
\caption{{\bf EEPA$^+$ Adaptive Sampling}}
\begin{algorithmic}[2]
\STATE $B:$Total Budget, $q:$Batch Size, $R:$Remained Budget, $S:$ Selected Subset 
\STATE $\mathcal{D}=\{\mathbf{x}_1, \ldots, \mathbf{x}_{n}\}$,  $\mathcal{F}= \{f(\mathbf{x}_i)|\mathbf{x}_i\in \mathcal{D}\}$ 
\STATE $R=B$, $S=\varnothing$
\WHILE {$R>0$}
 \STATE \emph{Modeling}: Fit a surrogate model $\hat{f}$ (e.g., GP) on $(\mathcal{D},\mathcal{F})$
  \STATE \emph{Coordinate selection}: Randomly select a decreasing subset of dimensions
  \STATE \emph{Candidate generation}: Perturb the selected coordinates of  $\mathbf{x}^o=\argmin_{\mathbf{x} \in \mathcal{D}}{f(\mathbf{x})}$ to obtain $\mathcal{X}$
  \STATE \emph{Pareto frontier}: Construct the non-dominated set $\mathcal{C}$ over the $\mathcal{X}$,  using $\hat{f}$ and $\Delta(\mathbf{x})_{x \in \mathcal{X}}$
     \STATE \emph{Sampling}: Select a subset of $q$ candidate points $S$ from $\mathcal{C}$
 \STATE \emph{Evaluation}: $\mathcal{F}_S=\{f(\mathbf{x}_i)|\mathbf{x}_i\in S\}$
  \STATE $\mathcal{D}= \mathcal{D}\cup \mathcal{S}$; $\mathcal{F}=\mathcal{F}\cup \mathcal{F}_S$
  
  
\ENDWHILE
\STATE Return $\mathbf{x}^{o}\in\argmin_{\mathbf{x}\in \mathcal{D}} f(\mathbf{x})$
\end{algorithmic}\label{alg:EEPA}
\end{algorithm}
Algorithm~\ref{alg:EEPA} represents the main steps of our proposed \textbf{EEPA$^+$}. Similar to most BBO strategies, \textbf{EEPA$^+$} starts with
fitting a surrogate model (e.g., GP) on the already evaluated data points $\mathcal{D}$ followed by a candidate generation (discretization) step to create an effective representation of the solution space. In this step, we dynamically select a subset of coordinates of the best-known solution and perturb to generate new candidate points similar to \textbf{DYCORS}. Next, we construct a Pareto frontier on the candidate points based on a \emph{distance} $\Delta(\mathbf{x})_{x \in \mathcal{X}}=\min_{x \in \mathcal{X}}{||\mathbf{x}-x||}$ and \emph{estimated function value $\hat{f}$} criterion and obtain a non-dominated set of points. To select a batch of points from the non-dominated set, \textbf{EEPA$^+$}
includes the point with minimum $\hat{f}(\mathbf{x})$
and eliminates the close points by applying a \emph{maximin} measure to identify the most diverse candidates.
In general, \textbf{EEPA$^+$} is able to adapt itself to both large and small batch sizes. When the non-dominated set is smaller than the batch size $|\mathcal{C}| \leq q$, it selects all available non-dominated points as the selected subset $S$ and saves the remaining evaluation budget $B-|S|$ for future iterations.



\begin{wrapfigure}{r}{0.45\textwidth} 
\centering
\vspace{-2mm}
\subfigure{\includegraphics[width=0.45\textwidth]{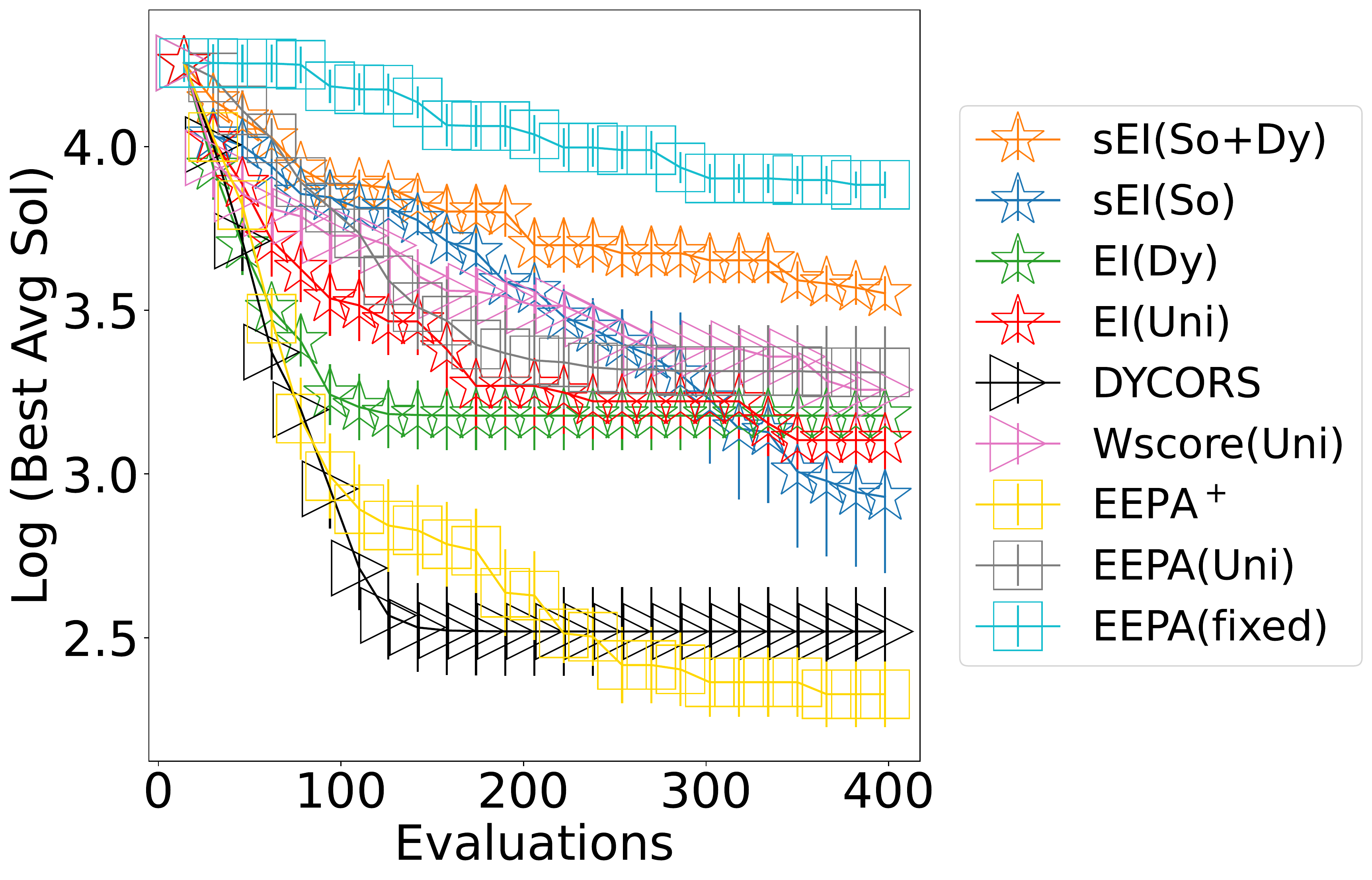}}
\vspace{-8mm}
\caption{Impact of discretization on the acquisition performance for Rastrigin function}
\label{exp:discretization}
\vspace{-4mm}
\end{wrapfigure}

Figure~\ref{exp:discretization} shows the performance of Expected Improvement (\textbf{EI}), \textbf{EEPA}, and Weighted Score (\textbf{Wscore}) acquisition functions utilizing different discretization approaches (uniform random (\textbf{Uni}, Sobol (\textbf{So}) \shortcite{owen1998scrambling}, and Dynamic coordinate search (\textbf{Dy})) for $6$-dimensional Rastrigin test problem. 
We can observe that the space discretization schemes have a significant impact on both the convergence pattern and the quality of the final solution obtained by distance-based acquisition functions (\textbf{EEPA} and \textbf{Wscore}). We refer to \textbf{Wscore} acquisition with \textbf{Dy} discretization as \textbf{DYCORS} in Figure~\ref{exp:discretization} and throughout the \S~\ref{sec:exp} and \S~\ref{sec:res}. On the other hand, the most widely-used variance-based acquisition function
\textbf{EI} 
is less sensitive to the choice of candidate generation strategy.
Note that the dynamic discretization approaches may not be necessarily as helpful for variance-based approaches (e.g., \textbf{EI}).
One justification for this observation is random perturbation of coordinates of the current best-known solution denoted by \textbf{Dy} (proposed in \shortcite{regis2013combining}), leads to convergence to local optima with \textbf{EI} acquisition function (green line). As discussed earlier, the variance-based acquisition functions without explicit trade-off between exploration and exploitation have weaker exploration strength, and the reduced subspace search of \textbf{Dy} restricts the exploration even further.
In this paper, we consider the best observed combination of acquisition functions and discretization. 

%% file: experiments.tex
\section{EXPERIMENTAL EVALUATION} \label{sec:exp}
In this section, we present a systematic comparison of the model-based adaptive sampling strategies in Black-Box Optimization literature as covered in \S~\ref{sec:sampling}. 
As discussed in \S~\ref{sec:sampling}, our baselines include sequential as well as the batch version of popular acquisitions for global optimization test problems.  Hence, we consider \textbf{sEI} and \textbf{qEI},  \textbf{sPI} and \textbf{qPI}, \textbf{sUCB} and \textbf{qUCB}, \textbf{sMES} and \textbf{qMES}, \textbf{qKG} denoting the sequential and batch versions of the discussed acquisition functions, respectively.
We mainly focus on the batch version for real-world application and hyperparameter tuning test problems due to their practicality.
Other considered baselines, \textbf{SOP}, \textbf{EEPA$^+$} and \textbf{DYCORS} are batch sampling techniques. In particular, for \textbf{DYCORS} we utilize the synchronous batch sampling proposed in \shortcite{eriksson2019pysot}. 
We consider Gaussian Processes (GP) with fixed hyperparameters as the surrogate model for all baselines across all considered test problems. More specifically, we utilize a 
GP, 
which uses a Constant mean and a Matérn Kernel with smoothness parameter $2.5$ and gamma length scale prior, $\Gamma \sim (\alpha=3,\beta=6)$ where $\alpha$ and $\beta$ represent concentration (shape) and rate parameters of the gamma distribution, respectively.

\subsection{Experiment Setup}
We first consider three commonly used synthetic optimization test problems \emph{Rosenbrock}, \emph{Rastrigin}, and \emph{Levy} \shortcite{simulationlib} for two low-dimensional ($d=6$) and high-dimensional ($d=30$) cases. 
For each synthetic test problem, we use Latin Hypercube Design (LHD) \shortcite{joseph2008orthogonal} with maximin criterion to generate 30 different initial sets of labeled points with the size $2*(d+1)$ and replicate our evaluation 30 times.
We consider a batch size of $|S|=4$ across all experiments for this group of test problems. When $d=6$, the budget and discretized solution space size are 400 and 1000, respectively. When $d=30$, we use a budget of 500, and discretized space size is 5000. 

We investigate two well-known benchmark problems in science and engineering domains,  \textbf{DNA binding} and \textbf{Airfoil design} optimization,
 as BBO baselines for real-world applications.
For \textbf{DNA binding}, \shortcite{barrera2016survey} utilizes protein-binding microarrays to assess the DNA binding activity of all feasible length-8 DNA sequences on 201 protein (human transcription factors) targets. Each transcription factor dataset could be viewed as a separate optimization problem. The objective is to identify DNA sequences with the maximum binding activity score.
We consider CRX and VSX1 protein binding problems similar to \shortcite{hashimoto2018derivative}.
Each letter of DNA code has been rewritten to a numeric value resulting in an $8$-dimensional array. The evaluation budget and size of the initial set are 1000 and 18, respectively. We consider a batch size of 100, similar to \shortcite{hashimoto2018derivative}, and replicate the results 10 times for each different initialization obtained from the LHD technique.

For \textbf{Airfoil design}, we use the UIUC database \shortcite{cfd3a4869ab4439490668729fc1ba7ef}, which provides the geometries of approximately 1,600 airfoil designs. The designs are then evaluated by the XFOIL 2-D aerodynamics simulator \shortcite{drela1989xfoil}. The objective is to maximize the lift divided by the drag fraction.
We consider a $10$-dimensional representation of airfoil for perturbation and candidate generation based on NURBS parameterization proposed in \shortcite{viswanath2011dimension}. 
The evaluation budget and size of the initial set are 300 and 1, respectively. In each iteration, we evaluate 30 airfoils ($|S|=30$), and replicate the results 10 times under the same initial point.

Finally, 
we consider three hyperparameter optimization problems for MNIST image classification with Logistic Regression (\textbf{MNIST-Logit}), Deep Neural network (\textbf{MNIST-DNN}), and Convolutional Neural network (\textbf{MNIST-CNN}). Across all problems, the evaluation budget and batch size equal 100 and 4, respectively, and we replicate 10 times for 10 different initial sets of size $d+1$ points where $d$ is the number of hyperparameters.
For \textbf{MNIST-Logit}, we consider  \emph{batch size}, \emph{learning rate}, \emph{L2 regularization}, and \emph{number of epochs} to be optimized. For \textbf{MNIST-DNN}, 
we consider \emph{batch size, learning rate, number of Dense units in the 1st layer, number of Dense units in the 2nd layer, dropout rate}, and \emph{number of epochs}. 
For \textbf{MNIST-CNN}, 
we consider \emph{batch size, learning rate, number of units in convolution layer, kernel size, number of dense units in the 1st layer, number of dense units in the 2nd layer, dropout rate}, and \emph{number of epochs}. Note that the ranges of hyperparameters considered are similar across all test problems. 
More specifically, batch size range of $(32,512)$, learning rate $(0,0.1)$, L2 regularization $(0,1)$, number of dense/convolution units $(32,512)$, kernel size of $(2,9)$, dropout rate $(0,0.9)$ and number of epochs $(0,10)$ has been considered in all problems.



%% file: Result.tex
\section{RESULTS}\label{sec:res}

\begin{figure*}[!htp]
\centering
\subfigure[RosenBrock $d=6$]{\label{exp:rosen-6d}\includegraphics[width=0.29\linewidth]{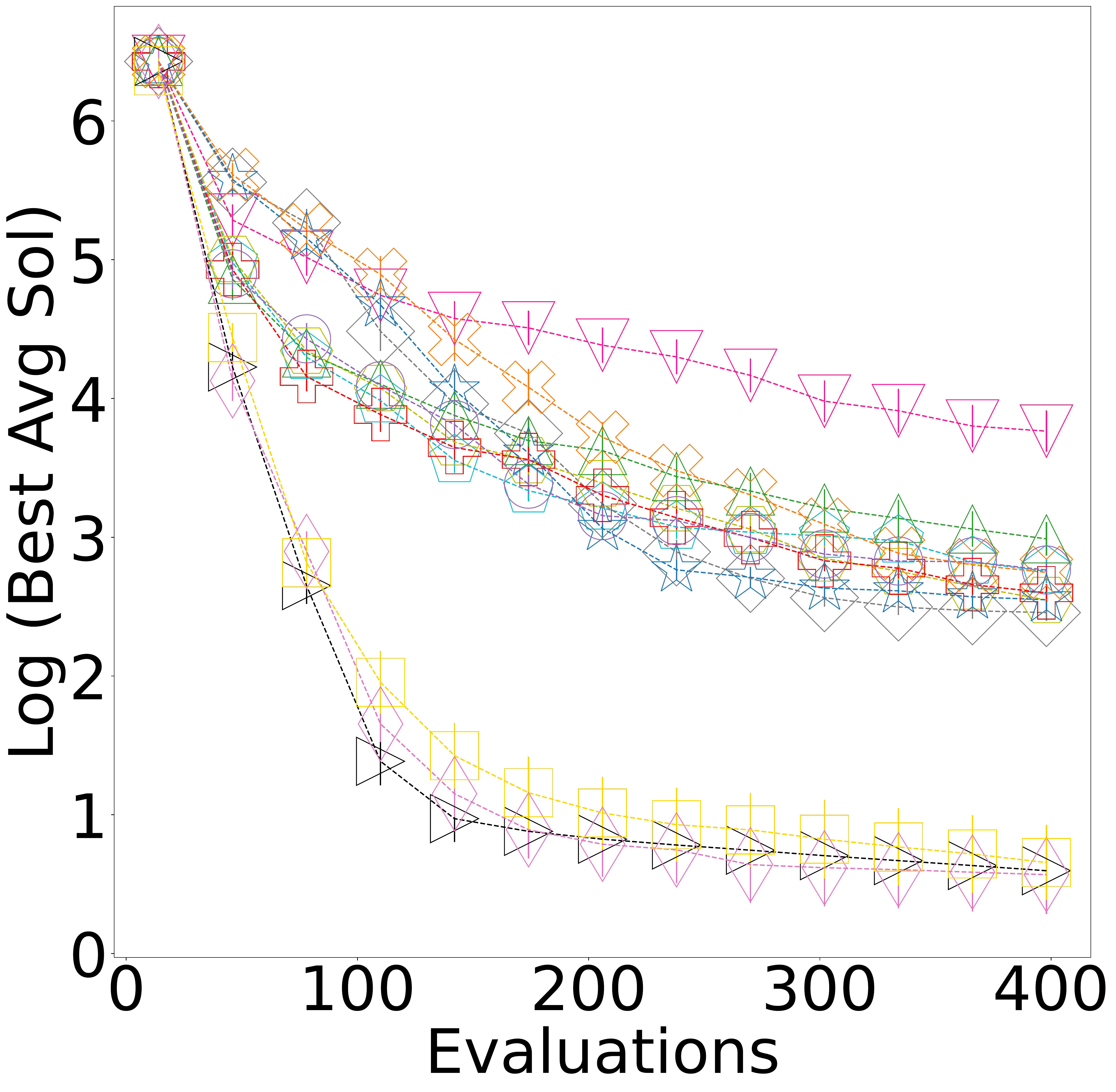}}
\subfigure[Rastrigin $d=6$]{\label{exp:rastrigin-6d}\includegraphics[width=0.29\linewidth]{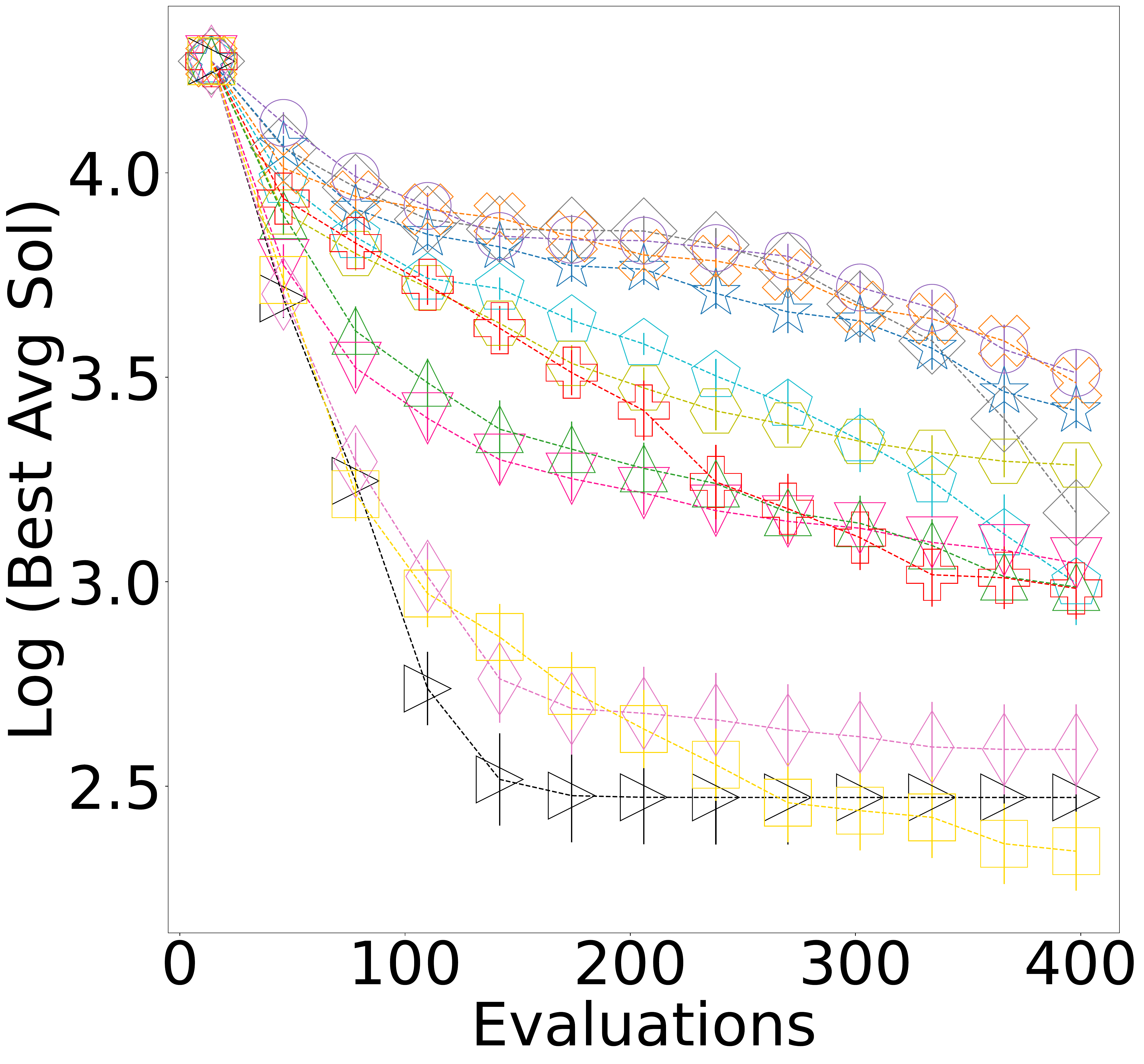}}
\subfigure[Levy $d=6$]{\label{exp:levy-6d}\includegraphics[width=0.29\linewidth]{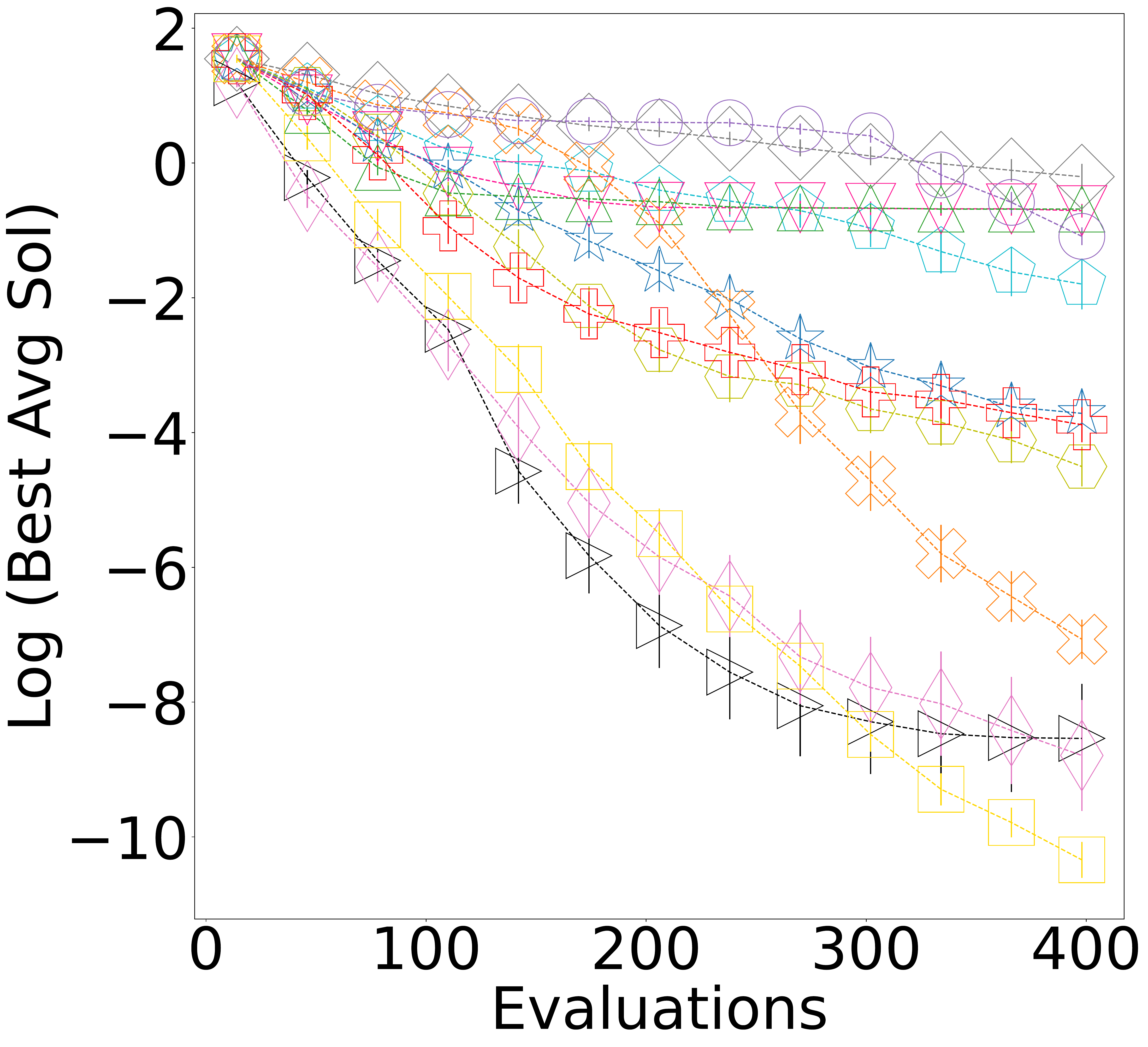}}
\vspace{-2mm}
\subfigure{\includegraphics[width=0.10\textwidth]{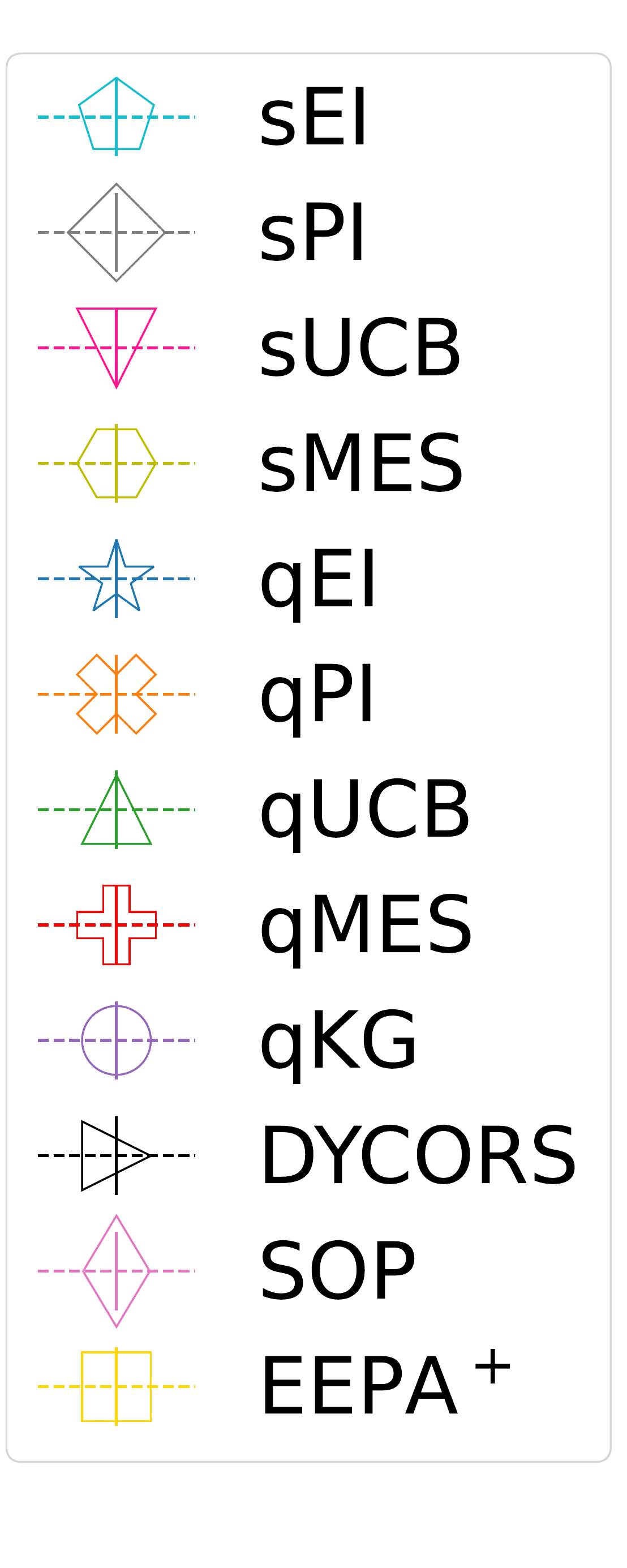}}
\subfigure[RosenBrock $d=30$]{\label{exp:rosen-30d}\includegraphics[width=0.29\linewidth]{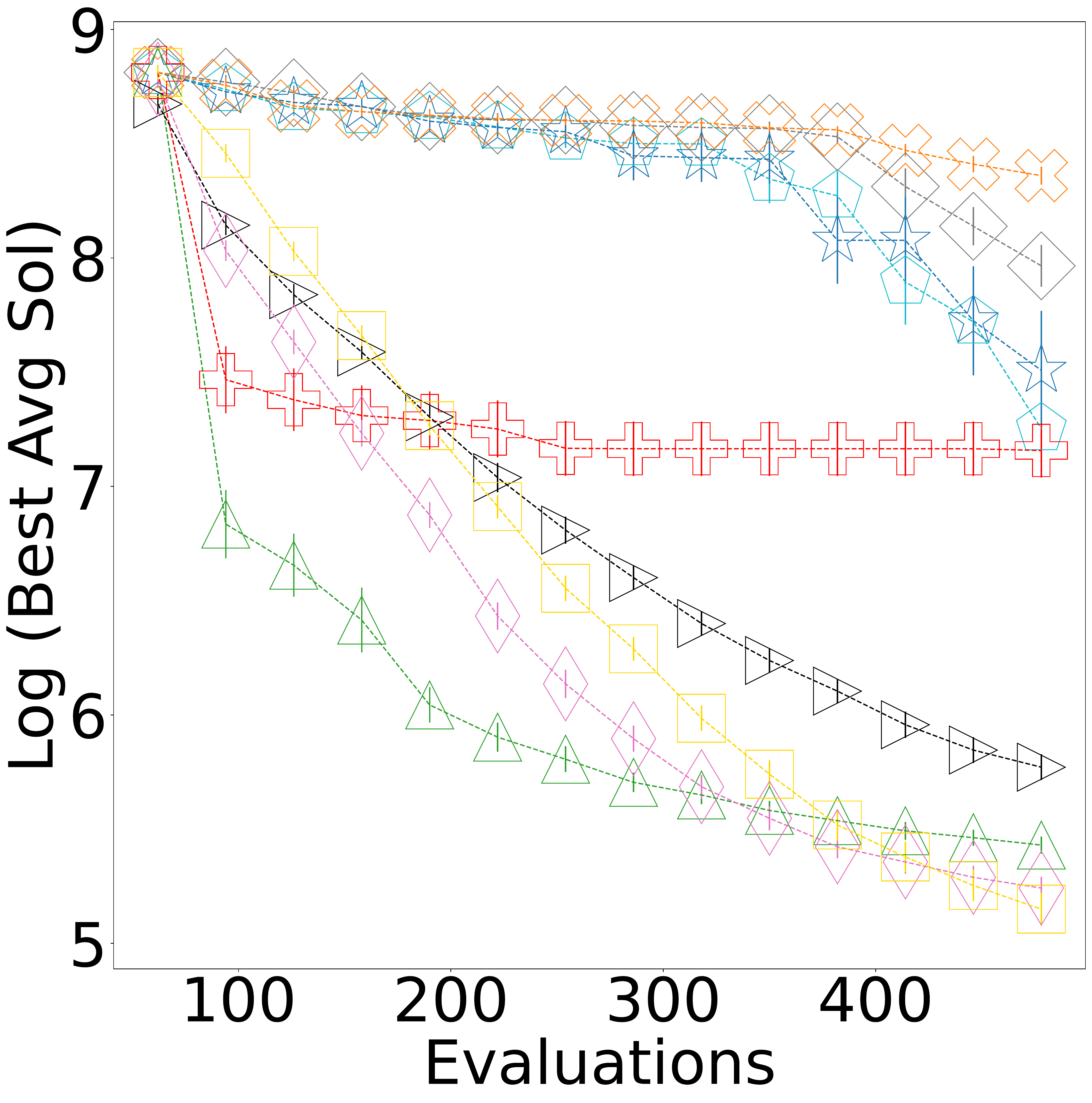}}
\vspace{-2mm}
\subfigure[Rastrigin $d=30$]{\label{exp:rastrigin-30d}\includegraphics[width=0.29\linewidth]{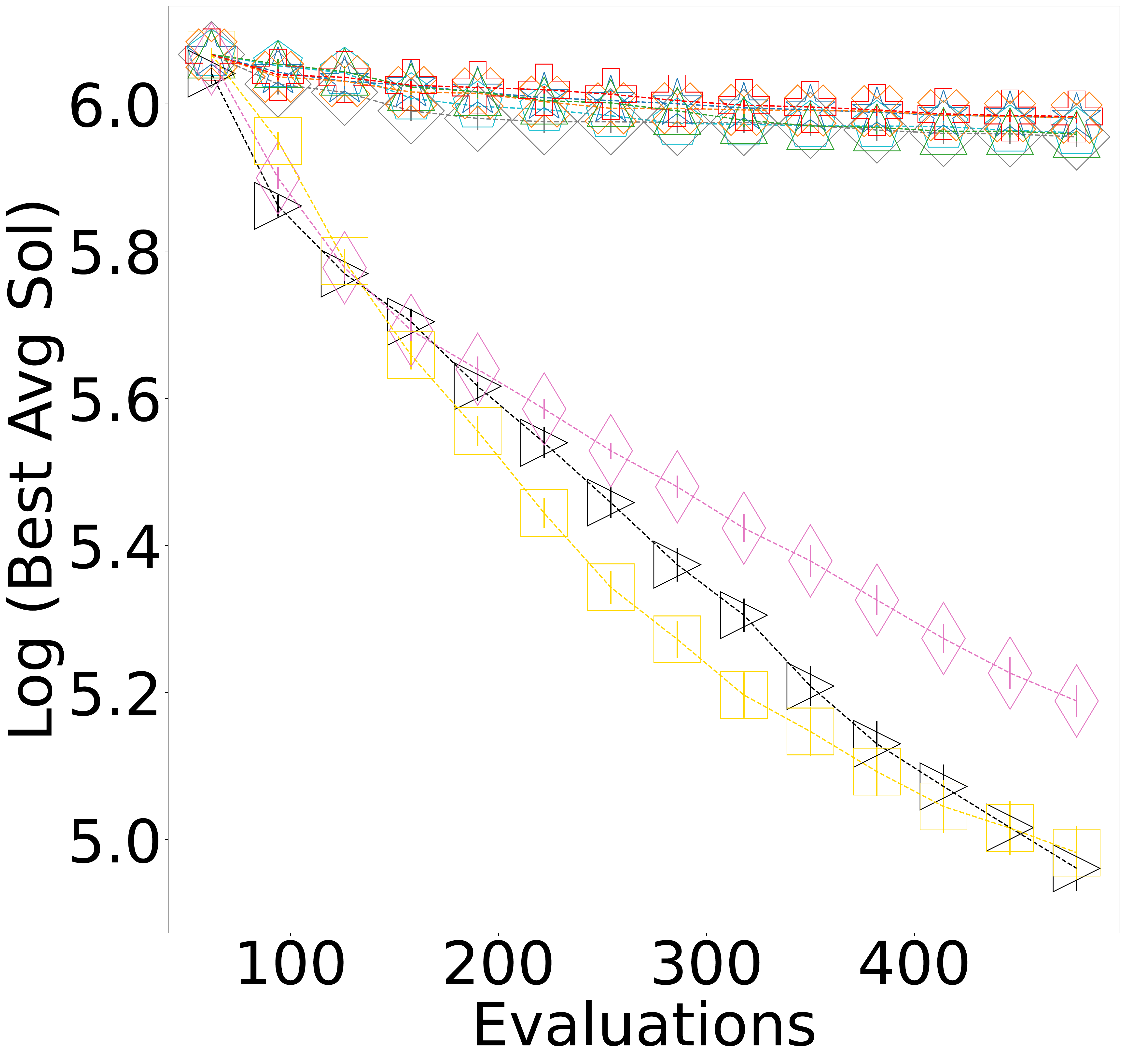}}
\vspace{-2mm}
\subfigure[Levy $d=30$]{\label{exp:levy-30d}\includegraphics[width=0.29\linewidth]{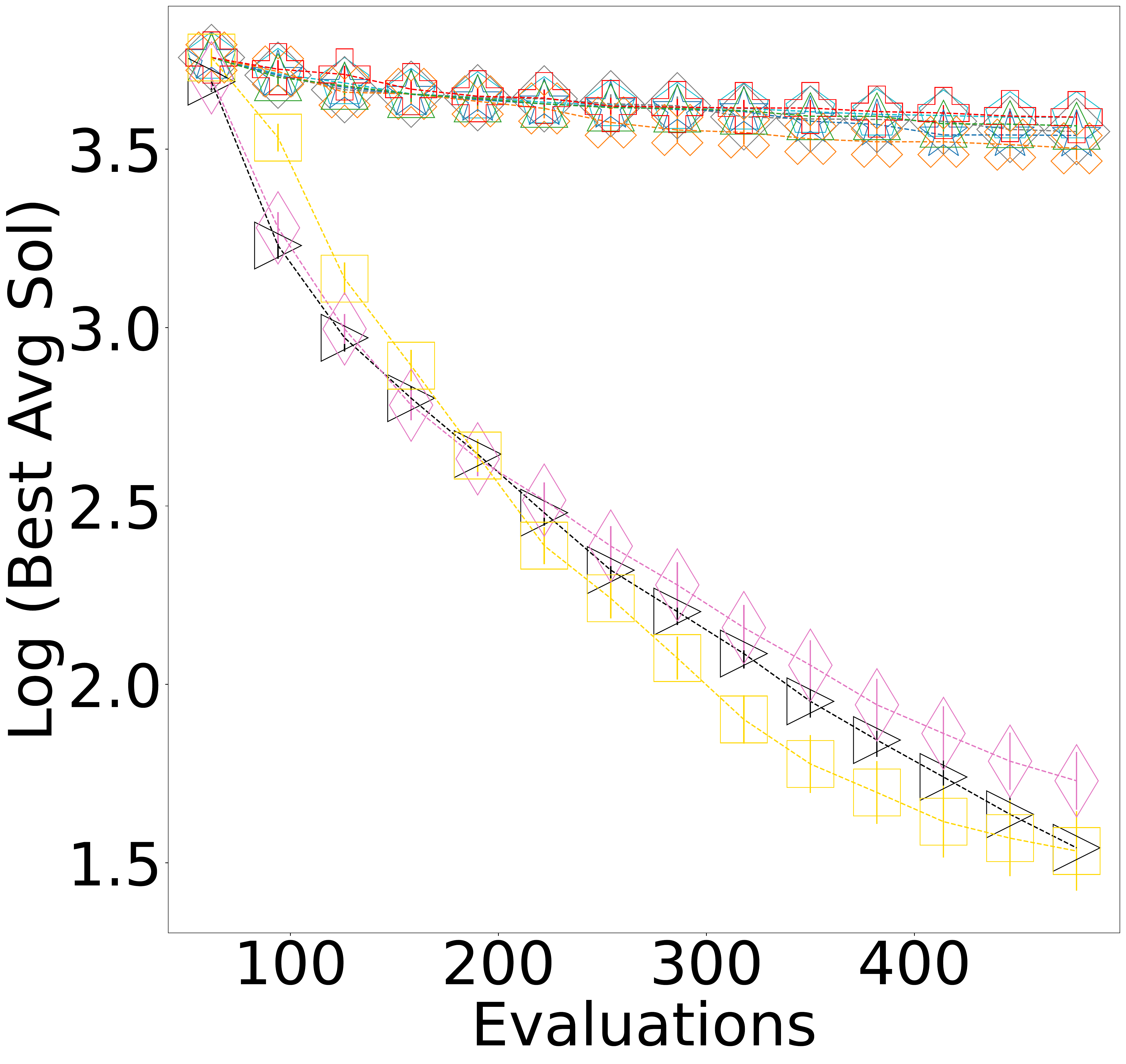}}
\vspace{-2mm}
\subfigure{\includegraphics[width=0.10\textwidth]{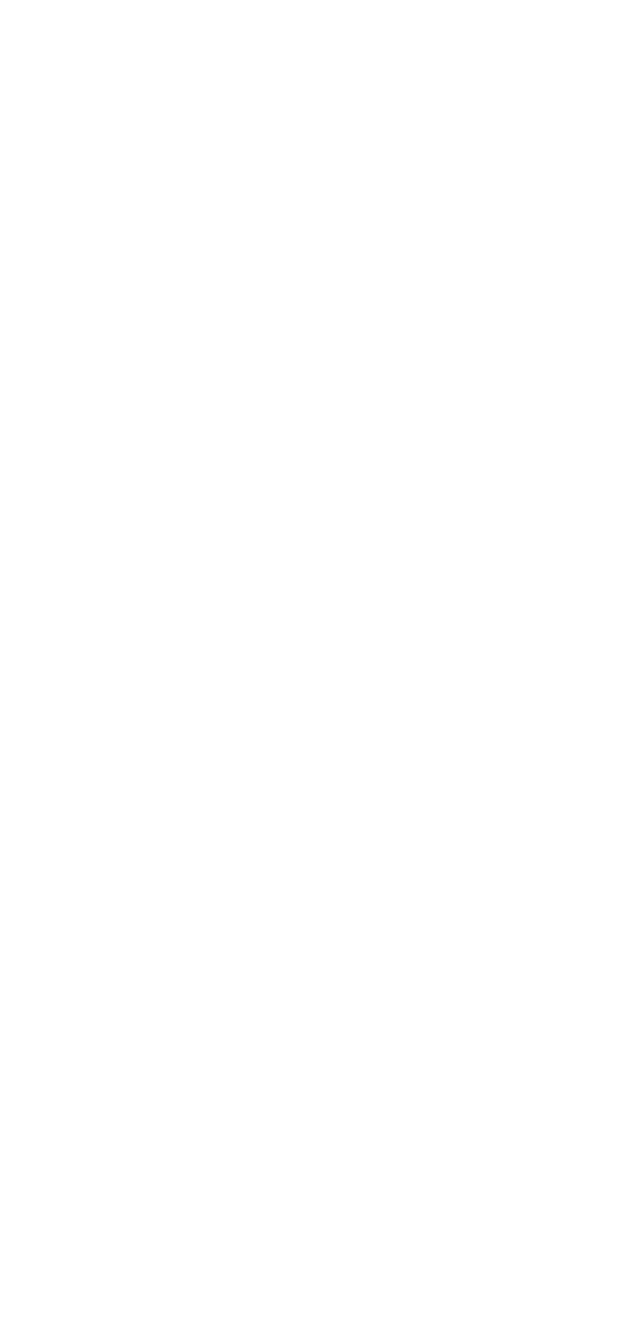}}
\vspace{1.5mm}
\caption{Global Optimization Test Problems} 
\label{exp:global}
\vspace{-8mm}
\end{figure*}

Figure~\ref{exp:global} demonstrates the performance of our considered baselines for the global optimization test problems plotted every 8 iterations.
It is worth noting that the performance of acquisition functions
diminishes as the dimension and complexity of problems increase. 
In general, distance-based acquisition functions (bf \textbf{EEPA$+$}, \textbf{SOP}, \textbf{DYCORS}) outperform variance-based acquisition functions for these problems, as a stronger global exploration is necessary for such complex functions. 
Note that for higher-dimensional problems, there are observable performance gaps between distance-based and variance-based baselines, making the choice of distance-based approaches even more practical.
In figure~\ref{exp:rosen-30d}, \textbf{qUCB} has a competitive performance with {\bf EEPA$^+$, SOP, DYCORS}. Although {\bf qUCB} is a variance-based acquisition function, it considers an explicit trade-off between exploration and exploitation criteria through the control parameter $\beta_{t}$.
Furthermore, 
{\bf qKG} is inapplicable for high-dimensional test cases due to the expensive computation.  

\begin{figure*}[!htp]
\centering
\subfigure[Binding to CRX Protein]{\label{exp:DNA1}\includegraphics[width=0.288\linewidth]{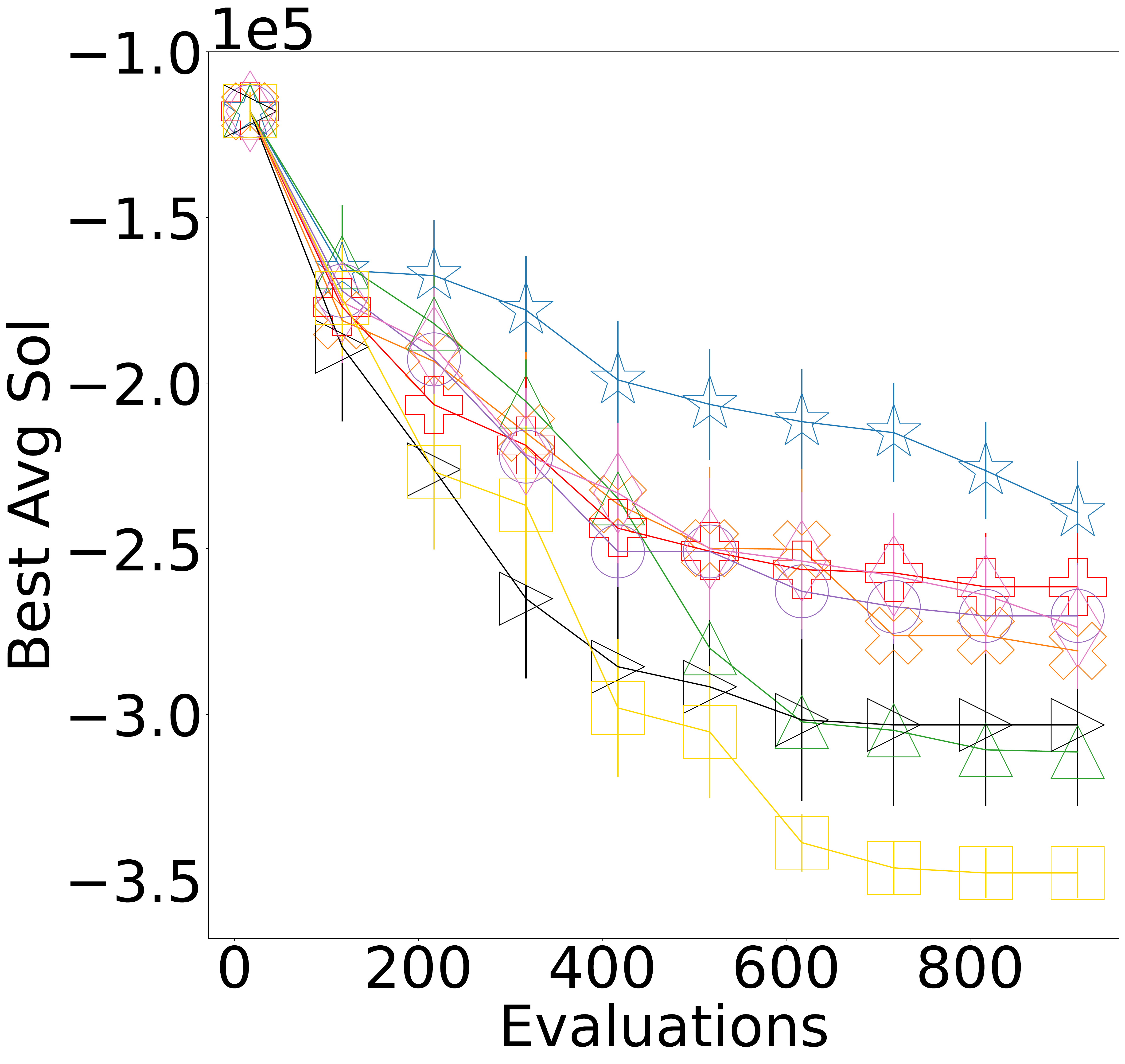}}
\subfigure[Binding to VSX1 Protein]{\label{exp:DNA2}\includegraphics[width=0.288\linewidth]{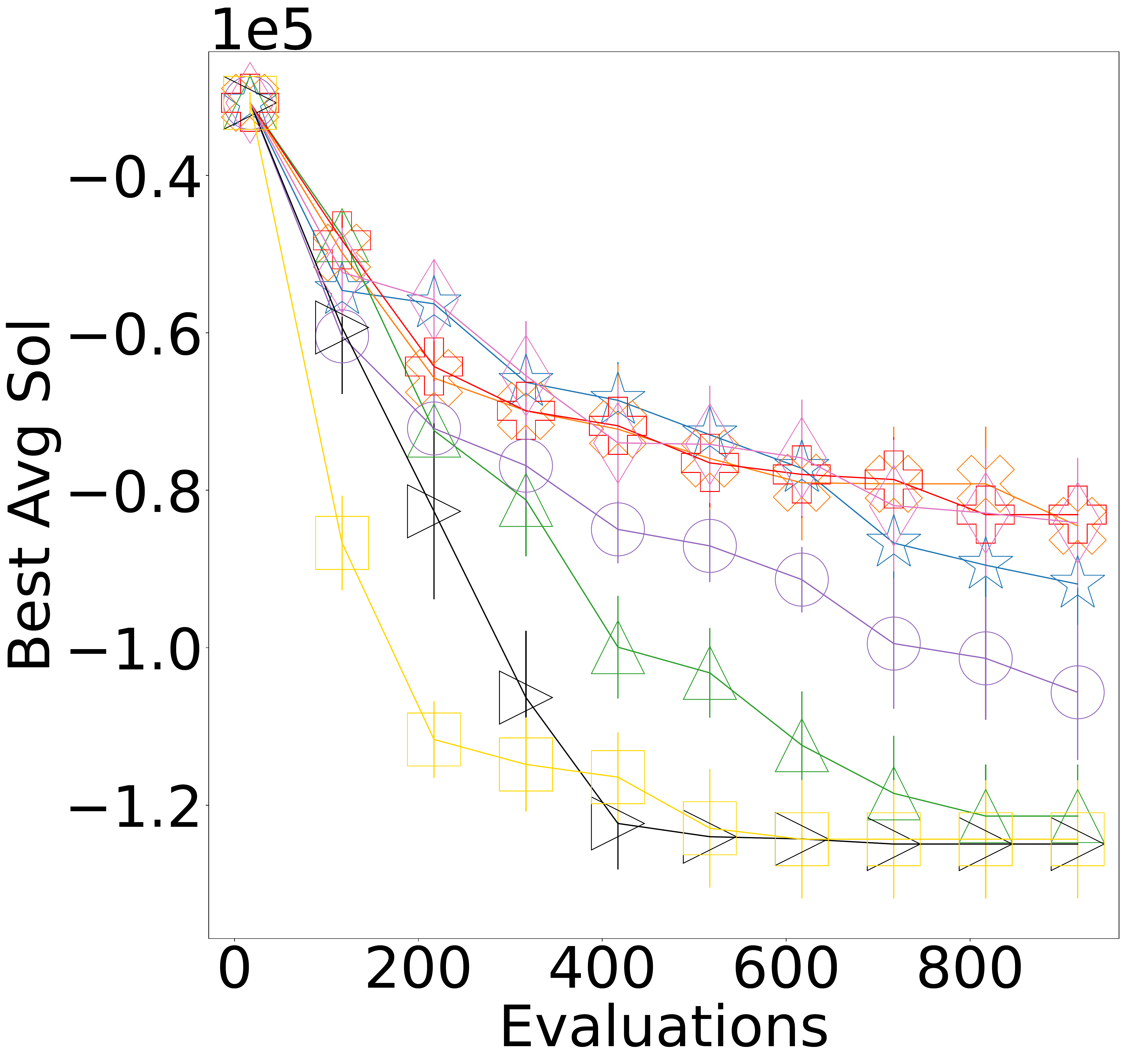}}
\subfigure[Airfoil Design]{\label{exp:airfoil}\includegraphics[width=0.288\linewidth]{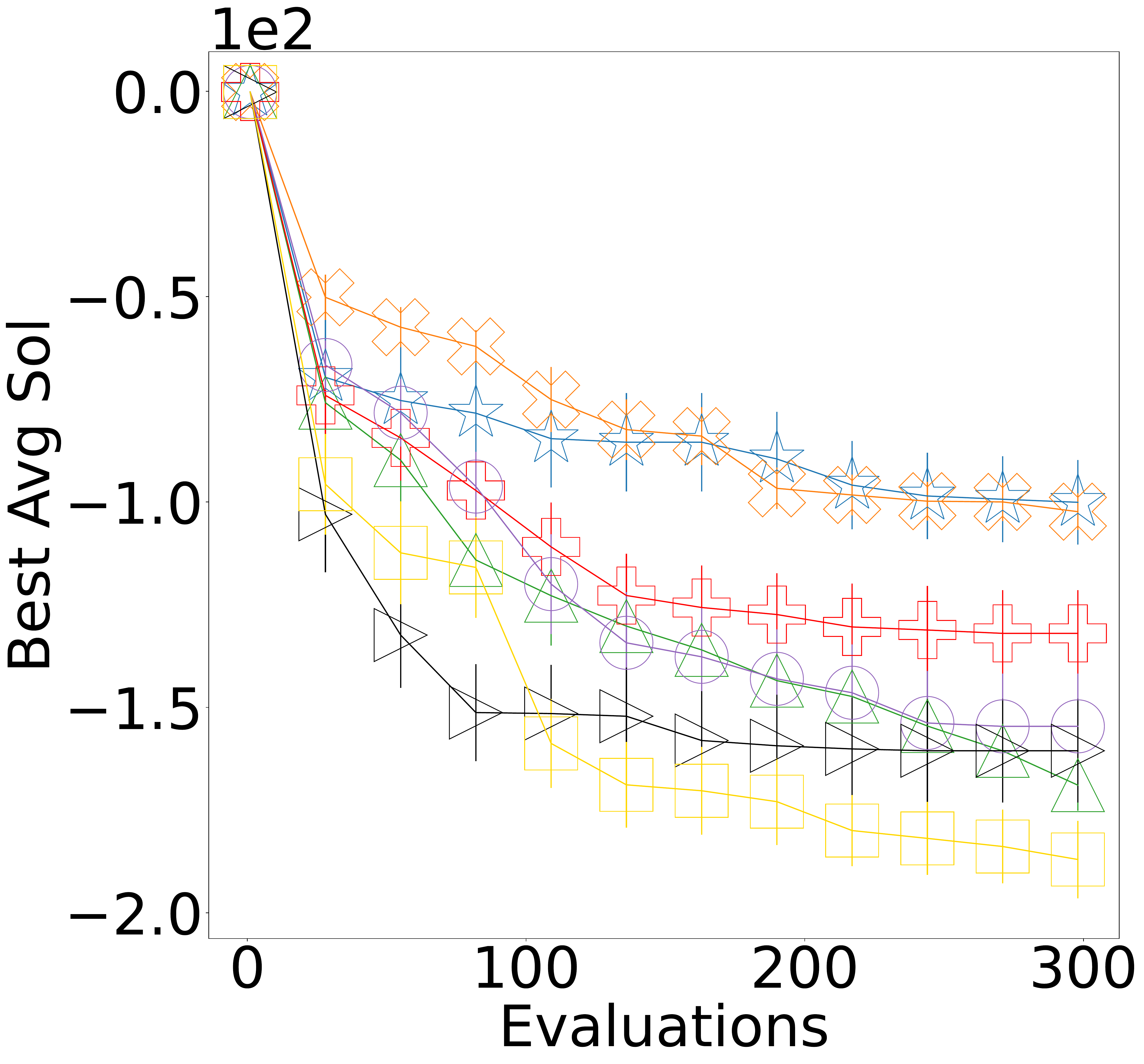}}
\subfigure{\includegraphics[width=0.11\textwidth]{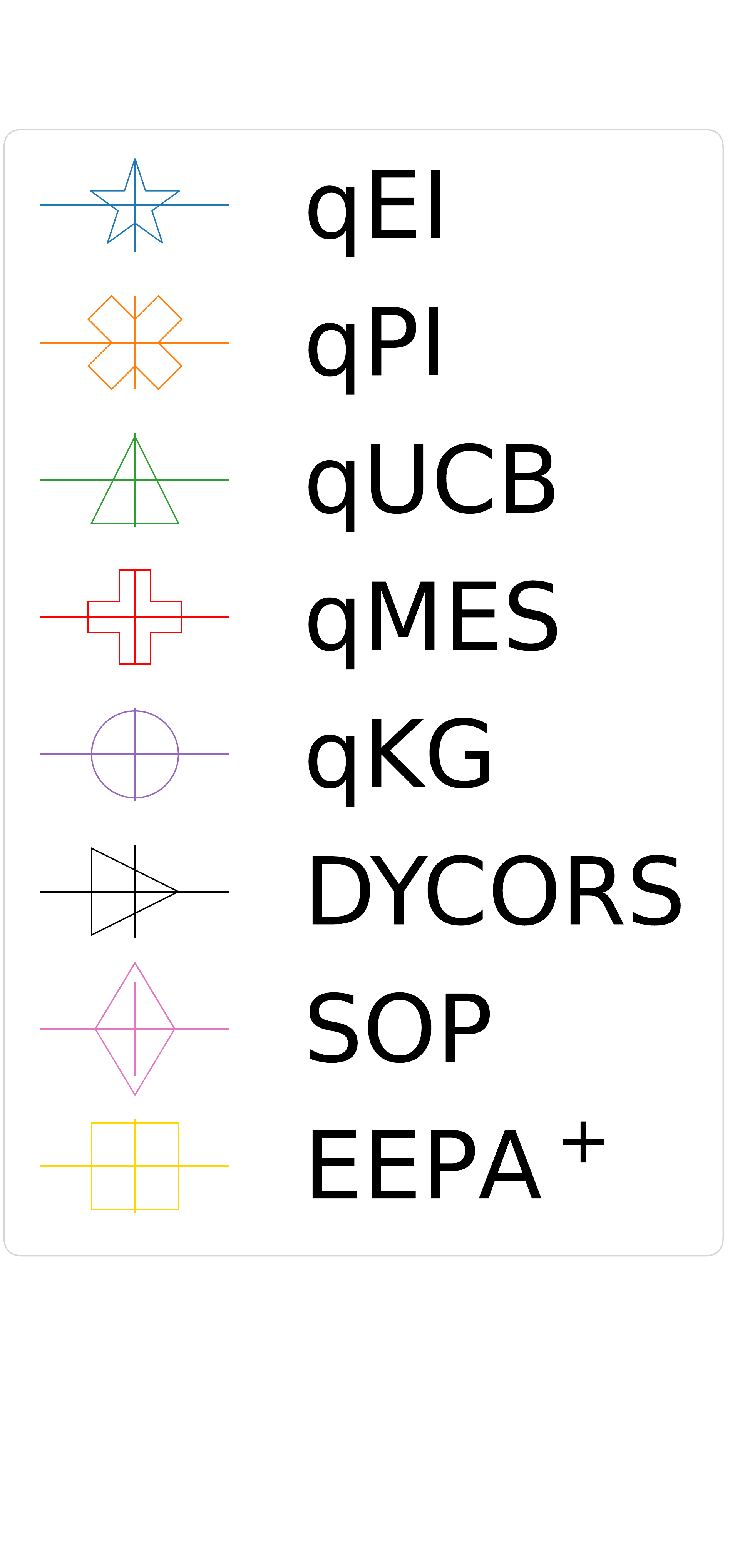}}
\vspace{-3mm}
\caption{Real world Optimization problems}
\label{exp:real}
\vspace{-7mm}
\end{figure*}

\begin{figure*}[!htp]
\centering
\subfigure[MNIST-Logit]{\label{exp:MNIST_Logit}\includegraphics[width=0.28\linewidth]{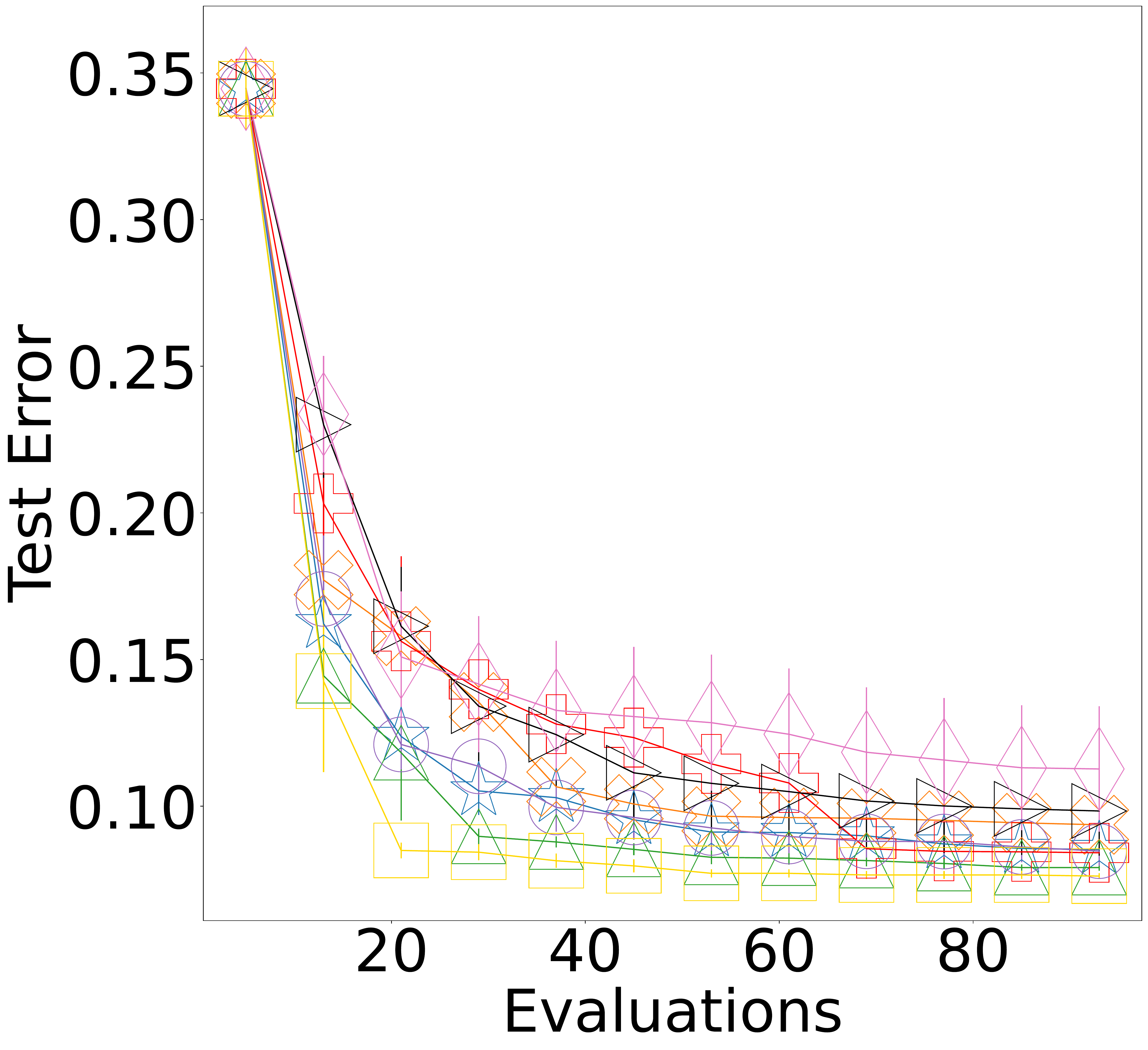}}
\subfigure[MNIST-DNN]{\label{exp:MNIST_DNN}\includegraphics[width=0.28\linewidth]{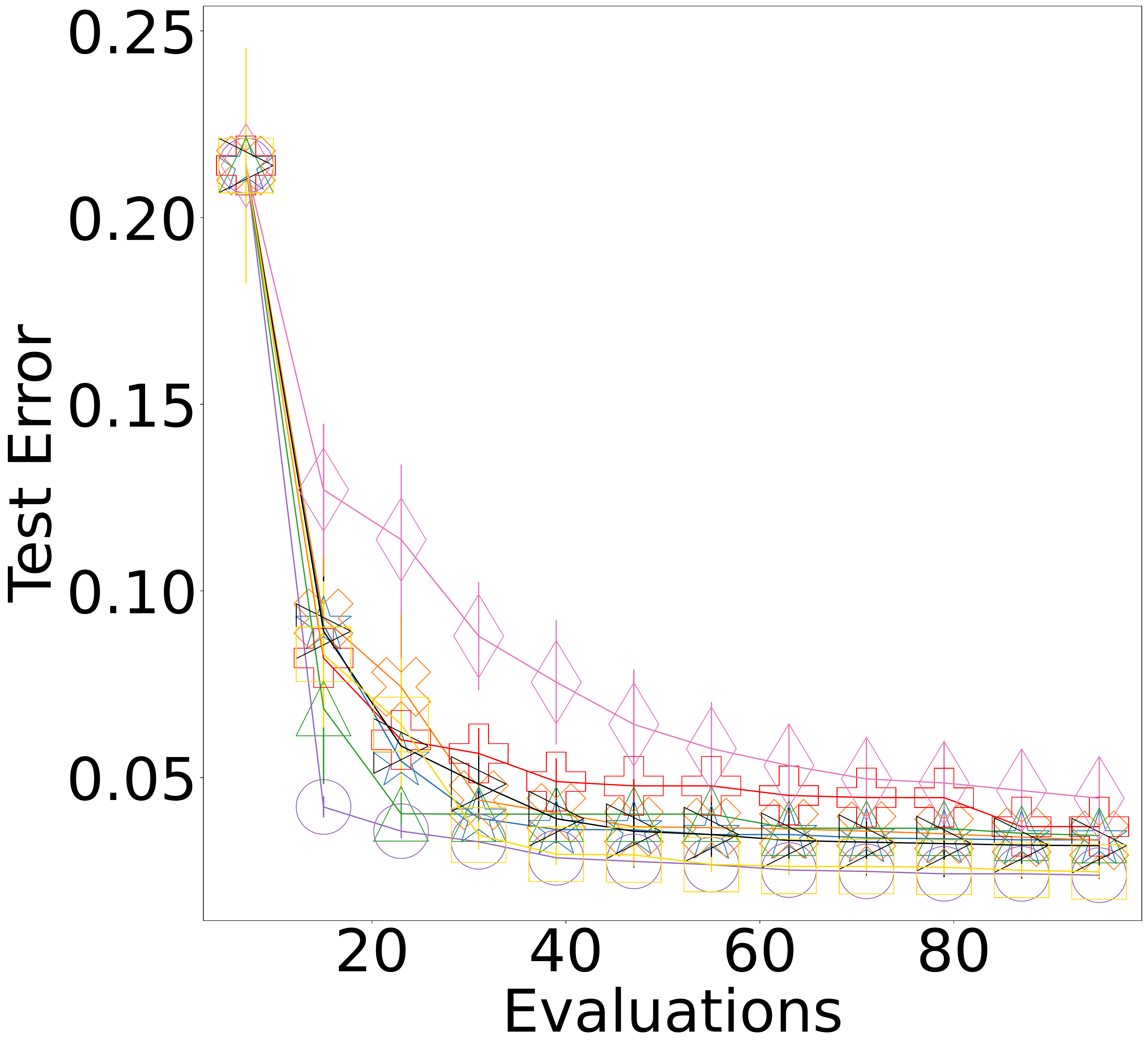}}
\subfigure[MNIST-CNN]{\label{exp:MNIST_CNN}\includegraphics[width=0.298\linewidth]{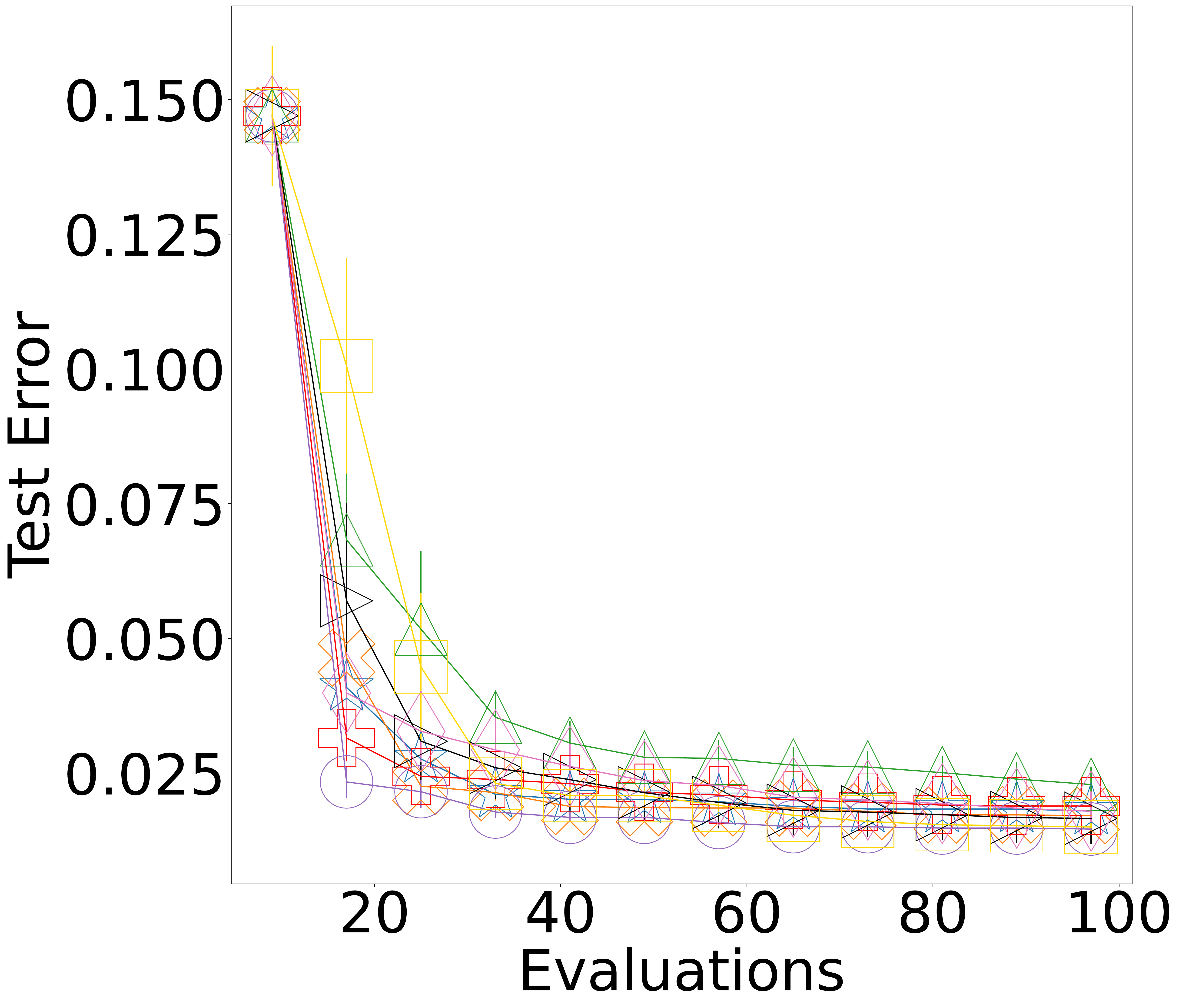}}
\subfigure{\includegraphics[width=0.11\textwidth]{fig/legend-DNA.pdf}}
\vspace{-3mm}
\caption{Hyperparameter Optimization}
\label{exp:hyperparam}
\end{figure*}

\begin{figure*}[!htp]
\centering
\subfigure[Synthetic Optimization ]{\label{exp:bx_synthetic}\includegraphics[width=0.305\linewidth]{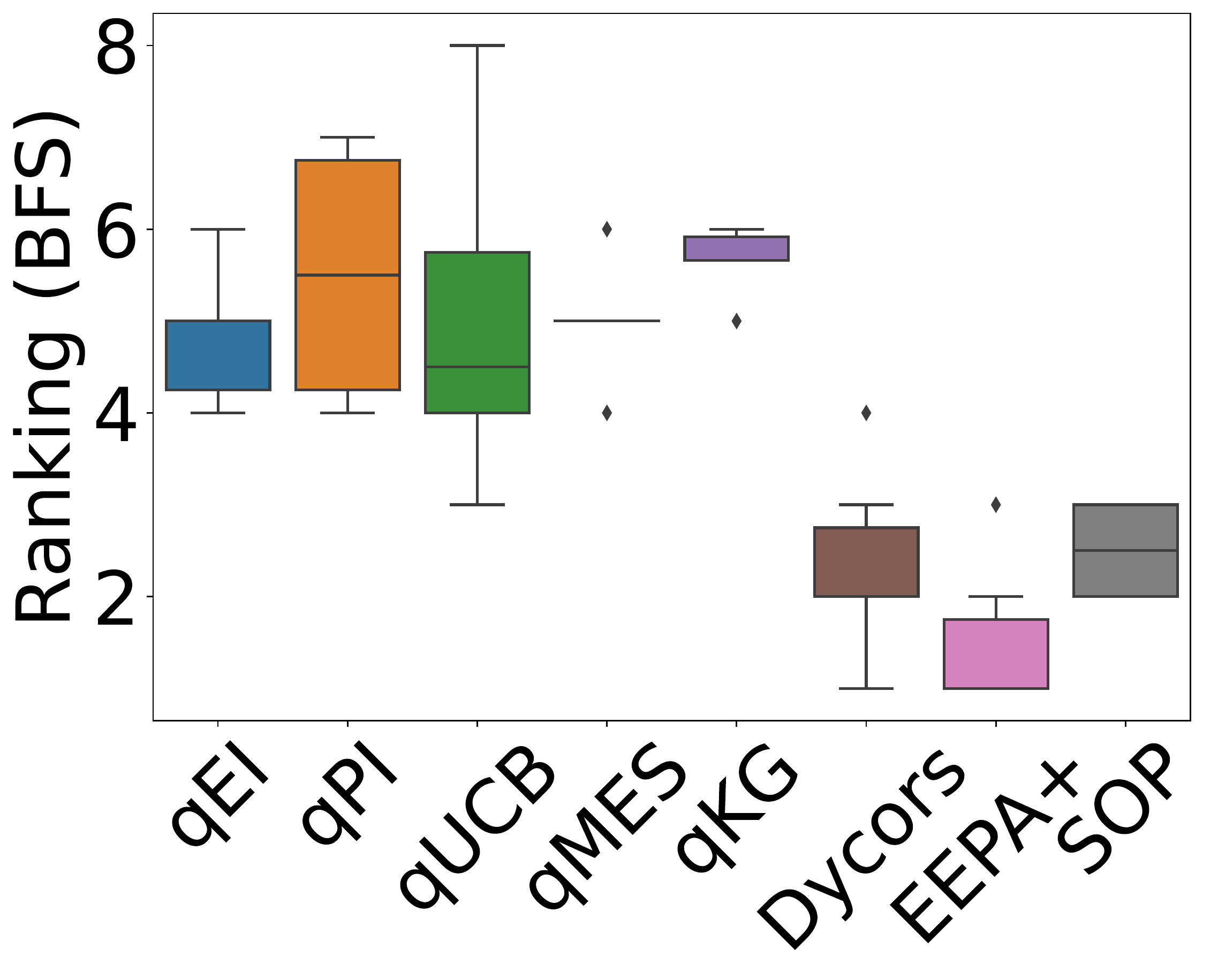}}
\subfigure[Real Optimization]{\label{exp:bx_real}\includegraphics[width=0.3\linewidth]{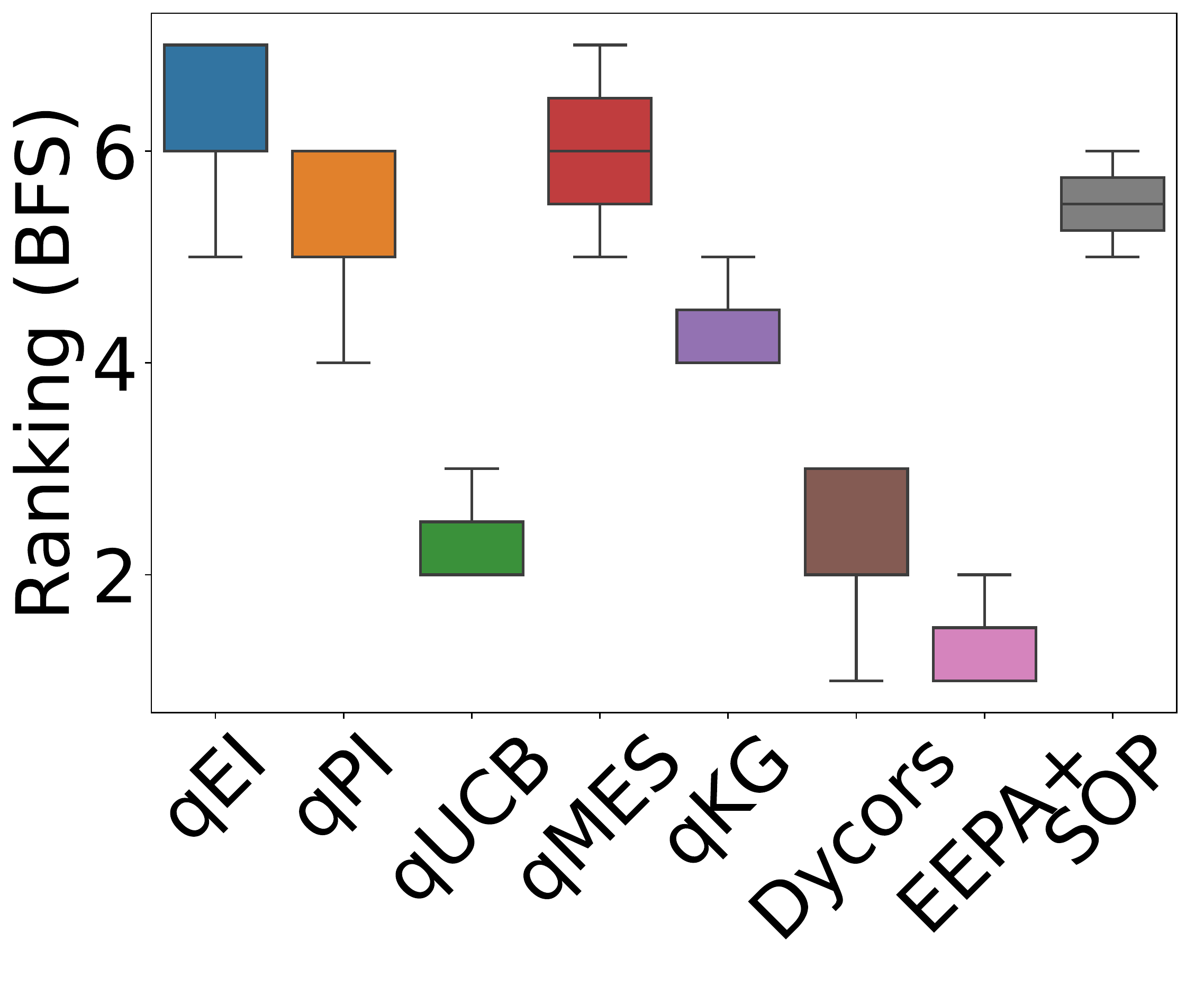}}
\subfigure[Hyperparameter Optimization ]{\label{exp:bx_hyper}\includegraphics[width=0.3\linewidth]{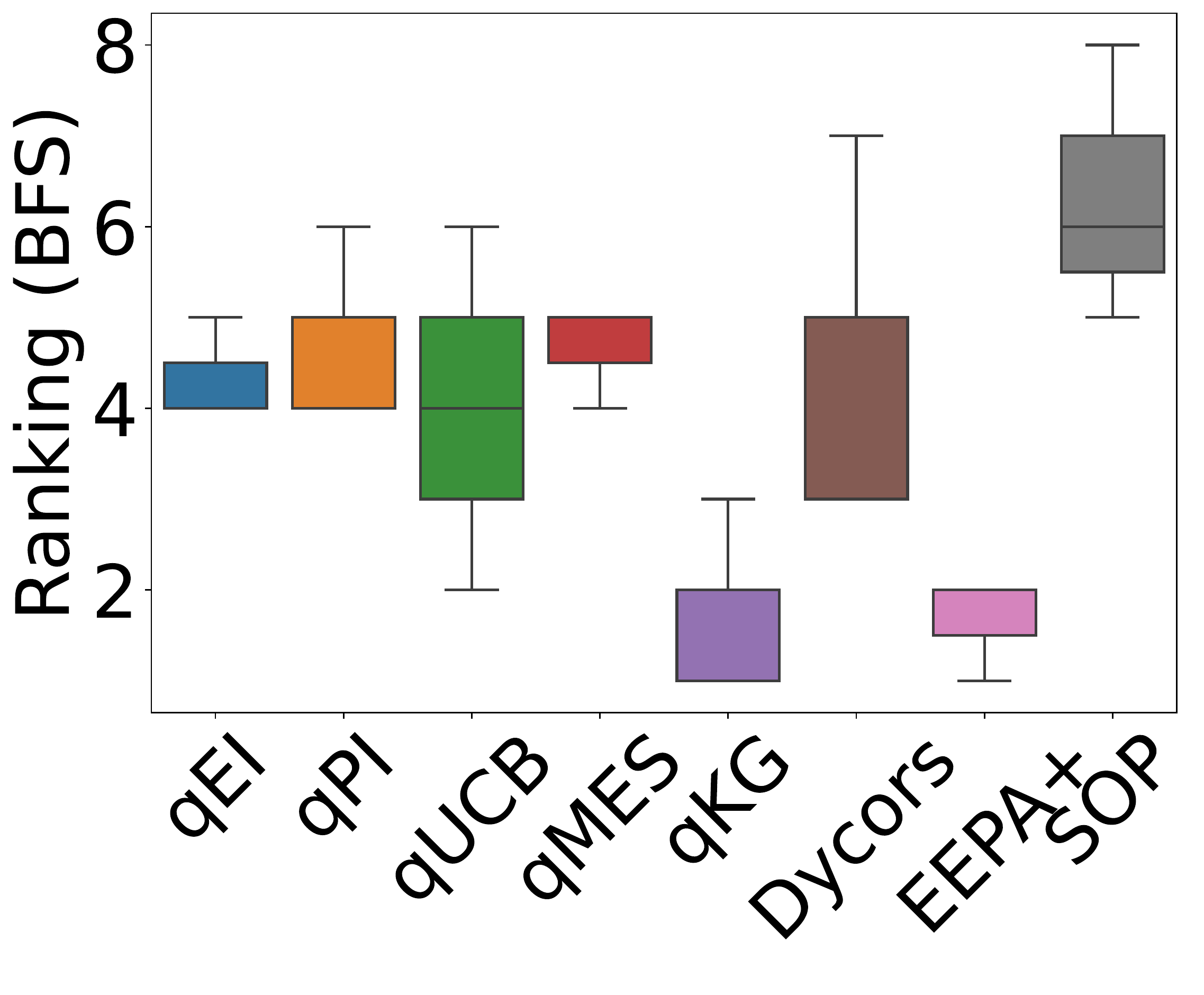}}
\caption{Ranking of BBO strategies across different problems} 
\label{exp:boxplots}
\vspace{-4mm}
\end{figure*}

In figure~\ref{exp:real}, we illustrate the  performance of the considered acquisition functions for real-world optimization test problems \textbf{DNA binding} and \textbf{Airfoil design}. Overall, we can observe a similar pattern as in global optimization test problems comparing distance-based versus variance-based acquisitions except for the {\bf SOP}. Although \textbf{SOP} has a promising performance for synthetic test problems, it lacks adequate flexibility to maintain its exploration power when batch size is large (e.g., 100). 
If the batch size is greater than the number of centers of \textbf{SOP}, it requires manipulation to adjust the size of the selected samples from each cluster, which adversely affects the performance of the algorithm. 
Moreover, {\bf SOP} considers multiple heuristic steps to filter out non-promising centers that need adjustments to be suitable for real-world applications. 
We can observe that distance-based approaches such as {\bf EEPA$^+$} and {\bf DYCORS} are still effective in this case. Note that {\bf qUCB}, has relatively promising convergence performance due to its adaptive consideration of exploration-exploitation trade-off parameter ($\beta_{t}$). 

Figure~\ref{exp:hyperparam} presents our empirical study on MNIST image classification task (per 8 iterations).
The results indicate that the variance-based acquisition functions are more competitive
for hyperparameter optimization. 
In particular, 
our proposed \textbf{EEPA$^+$} has competitive performance in all three hyperparameter tuning problems. {\bf qKG} also outperforms many other strategies for both \textbf{MNIST-DNN} and \textbf{MNIST-CNN}.
Unlike synthetic test cases, {\bf SOP} and {\bf DYCORS} are not as competitive for hyperparameter tuning. However, {\bf DYCORS} is in a reasonable range of performance.
The results indicate that our considered baselines are equally good for hyperparameter tuning problems since the response surface is easier to optimize as discussed in \S~\ref{sec:bbo}.

Finally, we assign a rank to each considered acquisition function in figure~\ref{exp:boxplots} based on the mean and variance of the final best-known solution obtained across all replications to summarize the results. 
We can confirm that
the variance-based acquisition functions are outperformed by the distance-based ones for complex optimization settings (both real applications and synthetic test problems), which makes distance-based approaches,  specifically \textbf{EEPA$^+$}, a more viable choice. For hyperparameter tuning problems, we can observe a more competitive performance of variance-based acquisition functions. In particular, \textbf{qKG} is competitive to {\bf EEPA$^+$}. However, our proposed {\bf EEPA$^+$} outperforms most of the acquisition functions across all test problems. {\bf EEPA$^+$} provides greater flexibility for the exploration-exploitation trade-off through the non-dominated set construction and the dynamic discretization scheme.

\section{FINAL REMARKS}\label{sec:final}
In this paper, 
we introduced new perspectives on evaluating 
the exploration component of acquisition functions and the discretization scheme in various settings. Our empirical and analytical analysis revealed that variance-based acquisitions (e.g., \textbf{qKG}) are more suitable for hyperparameter tuning settings. On the other hand, multi-criteria Pareto-based (e.g., \textbf{EEPA$^+$}), as well as weighted-score acquisition with dynamic discretization scheme (e.g., \textbf{DYCORS}), are more promising in finding the optimal regions of in high-dimensional complex optimization settings. 
Moreover, considering \emph{distance-based} exploration metric in acquisition functions yields promising results for engineering design problems. We also noticed that some of the existing adaptive acquisition functions lack the adequate flexibility of being conveniently adopted for different problems. \textbf{SOP} with a large batch size (which is necessary for certain problems) is not applicable and requires further adjustments. The well-known \textbf{qKG} is computationally inefficient with a large finite discretization scheme. 
Furthermore, an explicit trade-off for exploration and exploitation, similar to \textbf{qUCB}, can significantly improve the performance of variance-based approaches.
However, tuning the critical control parameter for weighted score acquisition functions including \textbf{qUCB} is a challenging task and requires extensive investigation. 

In this paper, we also investigated the impact of the discretization scheme on the performance of different acquisition functions. We discovered that discretization has a greater influence on distance-based acquisition functions than variance-based acquisition functions. Moreover, we noted that a dynamic discretization scheme is an effective approach for improving the performance of distance-based acquisition functions. As a result, we proposed an improved modification of an existing Pareto-based acquisition function. Our empirical analysis confirmed the effectiveness of \textbf{EEPA$^{+}$} in different settings.